\documentclass[11pt,openany]{book}
\usepackage{epsfig}

\newenvironment{proof}{{\it Proof:\/}}{$\Box$\vskip 0.08in}
\usepackage{amssymb}

 \newtheorem{theorem}{Theorem}[section]
 \newtheorem{lemma}[theorem]{Lemma}
 \newtheorem{corollary}[theorem]{Corollary}
 \newtheorem{remark}[theorem]{Remark}

\newtheorem{definition}[theorem]{Definition}

\newtheorem{conjecture}[theorem]{Conjecture}

\newtheorem{question}[theorem]{Question}
\newtheorem{example}[theorem]{Example}

\newtheorem{proposition}[theorem]{Proposition}

\newcommand{\pct}[1]{}
\newcommand{\Psfig}[1]{{\mbox{$\ \ $}}}
\newcommand{\bd}{\bf}
\newcommand{\R}{\mathbb R}

\begin{document}

\renewcommand{\thechapter}{\Roman{chapter}}
\thispagestyle{empty}


\
\vspace{0.5in}
 \begin{center}
 {\LARGE\bf KNOTS}\\
{\bf From combinatorics of knot diagrams to combinatorial topology
based on knots}
\end{center}

\vspace*{0.3in}

\centerline{Warszawa, November 30, 1984 -- Bethesda, October 31, 2004}
\vspace*{0.3in}

 \begin{center}
                      {\LARGE \bf J\'ozef H.~Przytycki}
\end{center}

\vspace*{0.3in}

\ \\
{\LARGE  LIST OF CHAPTERS}:\ \\
\ \\
{\LARGE \bf Chapter I: \ Preliminaries }\\
\ \\
{\LARGE \bf Chapter II:\ History of Knot Theory}\\
\ \\
{\LARGE \bf Chapter III:\ Conway type invariants }\\
\ \\
{\LARGE \bf Chapter IV:\  Goeritz and Seifert matrices}\\ \ \\
{\LARGE \bf Chapter V:\ Graphs and links}\\
{\bf e-print: http://arxiv.org/pdf/math.GT/0601227}\\ 
\ \\
{\LARGE \bf Chapter VI:\ Fox $n$-colorings, Rational moves, Lagrangian tangles
and Burnside groups}\ \\
\ \\
{\LARGE \bf Chapter VII:\ Symmetries of links}\ \\
\ \\
{\LARGE \bf Chapter VIII:\ Different links with the same
Jones type polynomials}\ \\
\ \\
{\LARGE \bf Chapter IX:\ Skein modules} \\
{\bf This e-print. Chapter IX starts at page 3}\\
\ \\
{\LARGE \bf Chapter X:\ Khovanov Homology: categorification of the Kauffman
bracket relation}\\
{\bf e-print: http://arxiv.org/pdf/math.GT/0512630 }\\ \ \\
{\LARGE \bf Appendix I.\ }\ \\
\ \\
{\LARGE \bf Appendix II.\ }\\ \ \\
{\LARGE \bf Appendix III.\ }\\
\
\newline

\ \\
{\LARGE \bf Introduction}\\
\ \\
This book is
about classical Knot Theory, that is, about
the position of a circle (a knot) or of a number of disjoint circles
(a link) in the space $R^3$ or in the sphere $S^3$.
We also venture into Knot Theory in general 3-dimensional
manifolds.

Lecture Notes on Knot Theory, published in Polish in 1995 \cite{P-18}, 
is the predecessor of this book\footnote{The
Polish edition was prepared for the ``Knot Theory" mini-semester
at the Stefan Banach Center, Warsaw, Poland, July-August, 1995.}.
A rough translation of the Notes (by J.~Wi\'sniewski) was
ready by the summer of 1995. It differed from the Polish edition
with the addition of
the full proof of Reidemeister's theorem. While I could not find
time to refine the translation and prepare the final manuscript,
I did add new material and rewrote existing
chapters. In this way I created a new book based on the Polish
Lecture Notes but expanded three-fold.
Only the first part of Chapter III (formerly Chapter II),
on Conway's algebras is essentially unchanged from the Polish book
and is based on preprints \cite{P-1}.

As to the origin of the Lecture Notes, I taught an advanced course
on the theory of 3-manifolds and Knot Theory at Warsaw University and it
was only natural to write down my talks (see Introduction to  (Polish)
Lecture Notes).
\\ \ \\
...\\
\ \  \ \ SEE Introduction before CHAPTER I.

\setcounter{chapter}{8}

\chapter{Skein modules}\label{IX}
\centerline{Bethesda, September 30, 2004}
\ \\ 

We  describe in this chapter the idea of building an
algebraic topology based on knots (or more generally on the position of
embedded objects). That is, our basic building blocks
are considered up to ambient isotopy (not homotopy or homology).
For example, one should start from knots in 3-manifolds, surfaces in
4-manifolds, etc. However our theory is, until now, developed only
in the case of links in 3-manifolds, with only a glance towards 
4-manifolds. The main object of the theory is a {\it skein module} and 
we devote this chapter mostly to the description of skein 
modules in 3-dimensional manifolds.
In this book we outline the theory of skein modules often giving
only ideas and outlines of proofs. The author is preparing a monograph
 devoted exclusively to skein modules and their ramifications
\cite{P-30}.


\section{ History of skein modules}\label{IX.1}

H. Poincar\'e, in his paper ``Analysis situs" (1895), abstractly defined 
homology groups starting from formal linear combinations of simplices,
choosing cycles and dividing them by relations coming from boundaries
\cite{Po}\footnote{Before Poincar\'e the only similar construction
was the formation of ``divisors" on an algebraic curve by Dedekind and
Weber \cite{D-W}, that is the idea of considering formal linear 
combinations of points on an algebraic curve, modulo relations yielded 
by rational functions on the curve.}.
 
The idea behind skein modules is to use links instead of cycles (in the
case of a 3-manifold).
More precisely we consider the free module
generated by links modulo properly chosen (local) skein relations.
 
Skein modules have their origin in the observation made by J.~W.~Alexander 
(\cite{Al-3}, 1928)\footnote {We can argue further that 
Alexander was motivated by the 
chromatic polynomial introduced in 1912 by George David Birkhoff 
\cite{Birk-1}. Compare also the letter of Alexander to O.~Veblen 
(\cite {A-V}, 1919) discussed in Chapter II.} that his polynomials of 
three links $L_+, L_-$ and $L_0$ in $S^3$ are linearly related 
(here $L_+,L_-$ and $L_0$ denote three
links which are identical except in a small ball as shown
in Figure 1.1). J.~H.~Conway rediscovered the Alexander observation and
normalized the Alexander polynomial so that it satisfies the skein relation
$$\Delta_{L_+}(z) - \Delta_{L_-}(z) = z\Delta_{L_0}(z)$$
 (\cite{Co-1}, 1969). In the late seventies, Conway advocated the idea of
considering  the free $Z[z]$-module over oriented
links in an oriented 3-manifold and
dividing it by the submodule generated by his skein relation
\cite{Co-2} (cited in \cite{Gi}) and \cite{Co-3}
(cited in \cite{Ka-1}). However, there is no published account of
the content of Conway's
talks except when $S^3$ or its submanifolds are analyzed.
The original name Conway used for this object was ``linear skein".\ \\

\centerline{\psfig{figure=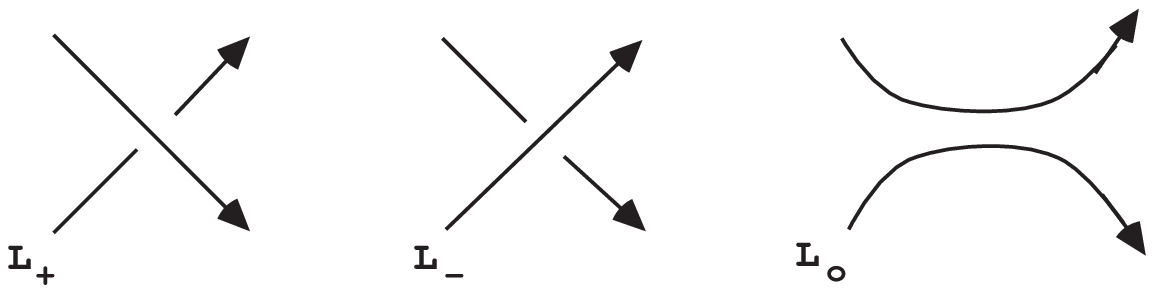}}
\begin{center}
Figure 1.1
\end{center}
 
Conway's idea was then pursued by Giller \cite{Gi}
(who computed the linear skein of a tangle), and Kauffman
\cite{Ka-1,Ka-8,Ka-3}, as well as Lickorish and Millett \cite{L-M-1} 
(for subspaces of $S^3$).

In graph theory, the idea of forming a ring of graphs and dividing 
it by an ideal generated by local relations was developed by W.Tutte in
his 1946 PhD thesis \cite{Tut-1}, but the relation to knot theory
was observed much later.

The theory of Hecke algebras, as introduced by N.~Iwahori (\cite{Iw},1964),
is closely connected to the theory of skein modules, however the relation
of Hecke algebras to knot theory was noticed by V.~Jones in 1984,
20 years after Iwahori's and Conway's work 
(it was crucial for Jones-Ocneanu construction of Markov traces).

The Temperley-Lieb algebra (\cite{T-L},1971) is related to the Kauffman
bracket skein module of the tangle, but any relation to knot theory was
again observed first by Jones in 1984.
 
At the time when I introduced skein modules, in April of 1987, I knew
the fundamental paper of Conway \cite{Co-1}, and \cite{Gi,Ka-8,L-M-1}
as well as \cite{Li-10}. However, the most stimulating paper
for me was one by J.~Hoste and M.~Kidwell \cite{Ho-K} about invariants 
of colored links, which I read in March of 1987. The goal of Hoste
and Kidwell was to find the common generalization of the multivariable
Alexander polynomial and the Jones-Conway (Homflypt) polynomial, 
and they succeeded, in the limited manner, requiring one color to represent 
a trivial component\footnote{At the first Cascade Mountain Seminar, 
in January of 1987, Jim Hoste explained to me his (then unfinished) work 
with Mark Kidwell and we realized that an analogous result for the Kauffman 
bracket is much easier to prove. We did not present, initially, our 
work in the language of skein modules \cite{H-P-1}.}.
I quickly realized that what  Hoste and Kidwell
actually computed were invariants of links in the solid torus,
$S^1\times D^2$, up to the Jones-Conway skein relation.
I wrote the introductory paper on skein modules in May of 1987
\cite{P-5}. Skein modules were also introduced independently
by V.~Turaev \cite{Tu-2}.

\section{Goal of skein module theory}\label{IX.2}
Our goal is to build an algebraic topology based on knots. We call the main 
object used in the theory {\it a skein module} and we associate it 
to any 3-dimensional manifold.  The essence is that
 skein modules are quotients of 
free modules over ambient isotopy classes of links in  3-manifolds 
by properly chosen local (skein) relations. 

These new objects are, after a slow start, intensively studied and 
their properties seem to be topologically very significant. 
In particular, one should look for their features that are similar 
to the Seifert-Van Kampen and Mayer-Vietoris theorems. 
Another interesting question concerns the relation between 
the skein modules of the base and the skein modules of the covering space, 
 for coverings and branched coverings. 
At present we can say something about the above question only in a very 
special situation and then results concern symmetries of links and 
3-manifolds. 
As in the case of homology, one should  try to understand  the free and 
torsion part of the module. In particular the torsion part of the module 
seems to reflect the geometry of the manifold (i.e. their incompressible 
surfaces). 

There is an ingenious argument that skein modules have, in general, a 
different nature than homotopy invariants. This was first observed 
by F.~Jaeger \cite{Ja-1} and explored in \cite{J-V-W}. 
Namely, computation of 
Jones type invariants is usually NP hard (and so, up to the famous conjecture 
cannot be performed in polynomial time \cite{G-J}) while computing, say,  
the Alexander polynomial or homology groups, can be accomplished
 in polynomial time. 

The most promising aspect of theory of skein modules is their 
interpretation as deformation of coordinate rings of character varieties.
In the case of Kauffman bracket skein module it is an $SL(2,C)$ 
character variety \cite{Bu-4,Bu-5,Bu-6,P-S-1,P-S-2} 
and for Homflypt skein module it is an
$SL(n,C)$ character variety \cite{Si-4,Si-5}.


\section{Skein modules of 3-manifolds; ideas and examples}\label{IX.3}
\markboth{\hfil{\sc Skein modules}\hfil}
{\hfil{\sc Skein modules: ideas and examples}\hfil}
 
Skein modules are quotients of free modules over ambient
isotopy classes of links in a 3-manifold by properly chosen local
(skein) relations.
The choice of relations is a delicate task as we
should take into account several factors:

\begin{enumerate}
\item [\textup{(i)}]  Is the module we obtain accessible (computable)?
\item [\textup{(ii)}]
How precise are our modules in distinguishing 3-manifolds and links
in them?
\item [\textup{(iii)}]
Does the module reflect topology/geometry of a 3-manifold
(e.g., surfaces in a manifold, geometric decomposition of a manifold)?
\item [\textup{(iv)}]
Does the module admit some additional structure (e.g., filtration,
gradation, multiplication, Hopf algebra structure)? Is it leading
to a Topological Quantum Field Theory (TQFT) by taking a finite
dimensional quotient?
\end{enumerate}
\noindent
From a practical point of view there is yet a fifth important factor
\begin{enumerate}
\item [\textup{(v)}]
The ``historical factor'' in the choice of (skein) relations:
the relations of links which were already studied (possibly for 
totally different reasons) will be compared with the new structures, 
just to see how they work in the new setup. 
For example, if we consider the Jones skein relation we
can be sure that even for $S^3$ we get a nontrivial result.
\end{enumerate}
 
The idea of the skein module should become more apparent after we consider
some examples.
 
\begin{example}[{\cite{P-5}}]\label{3.1}\
Let $M$ be an arbitrary 3-dimensional manifold and let  $R$
be a commutative ring with unit. Moreover, let ${\cal L}$ denote
the set of ambient isotopy classes of oriented links in $M$
and let $R\cal L$ denote the free $R$-module generated by $\cal L$.
In $R\cal L$ we consider the submodule  $S$ generated by
skein expressions of type $L_+-L_-$, where $L_+$ and $L_-$
are two oriented links in $M$ which are different in a small ball
where they are depicted in Fig.~3.1.
\ \\
\ \\
\centerline{\psfig{figure=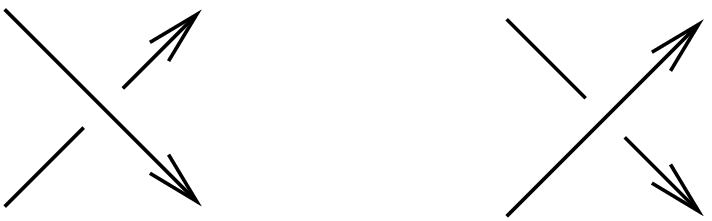}}
\begin{center}
Fig.~3.1
\end{center}
 
Let us note that, in order to define $S$, we do not assume 
that the manifold $M$
is oriented and there is no need to distinguish between $L_+$ and 
$L_-$ (sometimes it is not even possible to distinguish between them).
Now we define the skein module ${\cal S}(M;R,L_+-L_-)$,
which we denote simply by
${\cal S}_{\pm}(M)$, as the quotient $R{\cal L}/S$.
\end{example}
It is not hard to see that ${\cal S}_{\pm}(M)$ is a free module 
over the homotopy
classes of closed curves (links) in the manifold.
Let us have a closer look at this example.
The skein module ${\cal S}_{\pm}(M)$ admits a natural multiplication, 
defined as $L_1\cdot L_2=L_1 \sqcup  L_2$ (i.e.~the multiplication comes
from the disjoint sum of links).
Because of the skein relations we considered, the multiplication does not 
depend on the position of $L_2$ with respect to $L_1$.
The multiplication is associative and commutative but there is no
multiplicative identity. However, if we extend the set $\cal L$
by an empty link $\emptyset$ (and we call the resulting 
set ${\cal L}^{alg}$) then ${\cal S}(M)$ extends to 
the quotient ${\cal S}_{\pm}^{alg}(M)$
which is then a commutative algebra with an identity
(moreover the ring $R$ embeds into ${\cal S}_{\pm}^{alg}(M)$ 
by the mapping $r \mapsto r\cdot \emptyset$).
As an algebra,  ${\cal S}_{\pm}^{alg}(M)$ is isomorphic 
to the polynomial algebra
with coefficients in  $R$ and variables which are homotopy classes
of knots in $M$; in other words, the set of variables, $\hat{\pi}$,
is equal to the set of conjugacy classes of  the fundamental
group ${\pi}_1(M)$. Equivalently, the skein module
${\cal S}_{\pm}^{alg}(M)$ is an  $R$-algebra isomorphic to symmetric tensor
algebra over $R\hat{\pi}$ (which we usually denote by {\bd S}$R\hat{\pi}$).
Let us recall that the tensor algebra {\bd T}$R\hat \pi$ is the graded
sum $\bigoplus T^iR\hat \pi$
where $T^0R\hat \pi = R, T^1R\hat \pi = R\hat \pi, T^{i+1}R\hat \pi=
T^iR\hat \pi \otimes R\hat \pi$, and the symmetric algebra is the quotient
 ${\bd S}R\hat \pi = {\bd T} R\hat \pi /(a\otimes b -b\otimes a)$.
 
We consider this example in such detail since it turns out to be the
Vassiliev-Gusarov skein module of degree $0$.
The Vassiliev-Gusarov skein modules (which in fact are Hopf algebras)
give a good framework in the work on Vassiliev invariants of links
(in any 3-manifold); see Section 9 and \cite{P-9}.
 
Now, the above simple example can be generalized (or ``quantized'') to the
framed links, that is to the classes of annuli
embedded into a manifold (the central
curve of such an annulus determines an unframed link),
see Fig.~3.2.
\ \\
\ \\
\centerline{\psfig{figure=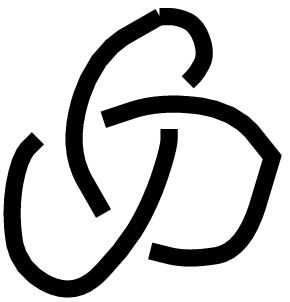}}
\begin{center}
Fig.~3.2
\end{center}

\begin{example}\label{3.2}.
Let  ${\cal L}^{fr}$ denote the set of ambient isotopy
classes of framed oriented links in an oriented 3-manifold  $M$.
Let  $R=Z[q^{\pm 1}]$. Let us consider the submodule $S^{fr}$
of $R{\cal L}^{fr}$ which is generated by skein expressions
presented in Fig.~3.3, that is,  $S^{fr}$
is generated by expressions of type $L_+-q^2L_-$ and $L^{(1)}-qL$,
where $L^{(1)}$ denotes a link obtained from  $L$ by one
positive twist of the framing of $L$.
 
\ \\
\ \\
\centerline{\psfig{figure=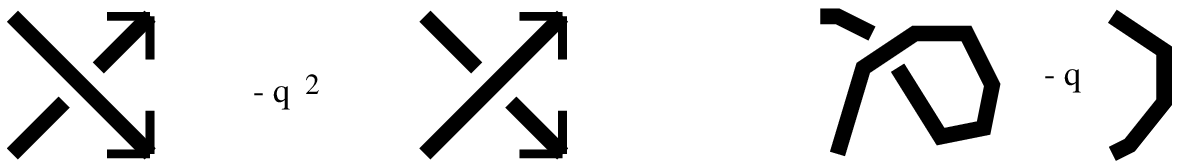}}
 
\begin{center}
Fig.~3.3
\end{center}

Now we define the skein module as the quotient ${\cal S}_{\pm}^{fr}(M)
= R{\cal L}^{fr}/ S^{fr}$.
\end{example}
 
The computation of ${\cal S}_{\pm}^{fr}(M)$ in general was an open problem 
for a long time.
The hard part is dealing with manifolds which contain a
non-separating torus (because such a torus is a source of a torsion in
the module). Initially I was only able to prove:
 
\begin{theorem}[\cite{P-19}]\label{1.3}
\begin{enumerate}
\item[(a)] If $M$ is an oriented 3-manifold without non-separating 2-spheres
and tori than ${\cal S}_{\pm}^{fr}(M)= {\cal S}(M,Z)\otimes Z[q^{\pm 1}]$.
\item[(b)] If $M$ contains a non-separating 2-sphere or torus then
${\cal S}_{\pm}^{fr}(M)$ has a torsion element:
\begin{enumerate} 
\item[(i)] If $L$ is a link in $M$ with the algebraic crossing number
with some 2-sphere in $M$ equal to $k$, $k\neq 0$, then $(q^{2k}-1)L =0$
in ${\cal S}_{\pm}^{fr}(M)$.
\item[(ii)] Let $L'$ be a link in $M$ with the algebraic crossing number 
with some torus in $M$ equal to $k$, $k\neq 0$. Let $L$ be a link obtained
by adding to $L'$ a noncontractible curve on the torus. Then $(q^{2k}-1)L =0$
in ${\cal S}_{\pm}^{fr}(M)$. 
\end{enumerate}
The link $L$ from (i) and (ii) is not equal to zero in the skein module 
because it is homotopy nontrivial and ${\cal S}_{\pm}^{fr}(M)$ reduces to
${\cal S}_{\pm}(M)$ when we set $q=1$.
\end{enumerate}
\end{theorem}

Uwe Kaiser completed description of ${\cal S}_{\pm}^{fr}(M)$ by 
proving that the only method of producing torsion is in fact 
described in (b)(i) and (b)(ii) \cite{Kai-4}.

In the example described above we used skein relations involving a 
crossing change. In the next example a smoothing of the crossing will 
be explored.

\begin{example}[{[H-P-1]}] \label{1.4}
Suppose that $M$ is an arbitrary manifold and $R$ is a commutative
ring with unit. Let $\cal L$ denote the set of isotopy classes of
oriented links in $M$ and $R\cal L$ denote a free $R$ module
over $\cal L$. Let $S_2$ be the submodule of $R\cal L$ generated by
skein expressions: \ \
\parbox{0.5cm}{\psfig{figure=L+maly.eps}} - 
\parbox{0.5cm}{\psfig{figure=L0maly.eps}}.\\
We define the second skein module ${\cal S}_2(M;R)$ to be the quotient
$R{\cal L} /S_2$.
\end{example}
 
\begin{lemma}[{\cite{H-P-1}}]\label{3.5}
Let $\phi: R{\cal L} \to RH_1(M,Z)$ be the epimorphism
of $R$ modules
such that ${\phi}(L)=|L|$ where $|L|$ is the class of $L$ in
$H_1(M,Z)$.\ Then $\phi$ can be factored through an isomorphism
$$\hat{\phi}: {\cal S}_2(M;R)\to RH_1(M,Z)$$
\end{lemma}
 
The idea of the proof is to show that: 

\begin{enumerate} 
\item 
[(i)] If two links, $L_1$ and $L_2$, represent the same element in 
$H_1(M,Z)$ then $L_1\sqcup -L_2$ is the boundary of an oriented 
surface contained in $M$, 
\item 
[(ii)] If $L_1\sqcup -L_2$ is the boundary of an oriented surface in 
$M$ then the link $L_2$ can be reached starting from $L_1$ via 
a sequence of elementary operations which are either modifications of 
\parbox{0.5cm}{\psfig{figure=L+maly.eps}} to 
\parbox{0.5cm}{\psfig{figure=L0maly.eps}} 
or their inversions. 
\end{enumerate} 

There is also a framed version of the second skein module  
and it can be called a $q$-deformation of the first homology group.
This version of the skein module contains torsion related to closed
surfaces in $M$ which do not separate $M$.
 
\begin{example}\label{3.6}
Let $M$ be an oriented  3-manifold and let
$R=Z[q^{\pm 1}]$. Let us denote by $S^{fr}$ the submodule of
$R{\cal L}^{fr}$ generated by skein expressions pictured in Fig.~3.4,
i.e.~$L_+-qL_0$ and $L^{(1)}-qL$, where $L^{(1)}$
denote a link obtained from  $L$ by twisting the framing of $L$ once 
in the positive direction.
 
\ \\
\ \\
\centerline{\psfig{figure=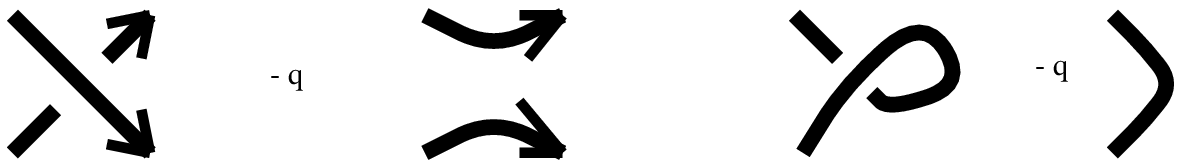}}
 
\begin{center}
Fig. 3.4
\end{center}
The quotient $R{\cal L}^{fr} /S^{fr}_2$ is called the second skein module and
it is denoted by ${\cal S}_2(M,q)$.
\end{example}
 
It is not hard to prove that ${\cal S}_2(M,q)$ is a free $R$ module
if $M$ is a rational homology sphere or its compact submanifold. In this
case ${\cal S}_2(M,q) = {\cal S}_2(M,Z)\otimes Z[q^{\pm 1}]$.
On the other hand, non-separating closed surfaces in $M$
yield torsion in the skein module.  We can compute the module
${\cal S}_2(M,q)$ for any 3-manifold.

\begin{theorem}[\cite{P-12}]\label{3.7}\ \\
$${\cal S}_2(M,q)=Z[q^{\pm 1}]T(H_1(M,Z))\oplus \bigoplus_
{\alpha \in H_1(M,Z)-T(H_1(M,Z))} Z[q^{\pm 1}]/(q^{2mul(\alpha)}-1),$$
where $T(H_1(M,Z))$ and the multiplicity of $\alpha$, $mul(\alpha)$, 
are defined as follows:
Let $\phi : H_1(M,Z)\times H_2(M,Z) \to Z$ be the bilinear form 
of the intersection of 1-cycles with 2-cycles on $M$
Then $\alpha \in T(H_1(M,Z))$ if and only if 
${\phi}_{\alpha}(\beta)=\phi(\alpha,\beta)=0$ for any $\beta$.
Otherwise $mul(\alpha)$ is defined as the positive generator of
$im{\phi}_{\alpha} (H_2(M,Z))$.
\end{theorem}

One should notice that for a closed 3-manifold $M$,
$T(H_1(M,Z))$ is the torsion part of $H_1(M,Z)$.

The second skein module should provide ``experimental data''
that is useful when computing harder skein modules.

We can consider even simpler skein module, without a skein relation,
but keeping the framing relation. We call such a skein module 
a framing skein module.

\begin{example}\label{3.8}
Let $M$ be an oriented  3-manifold and let
$R=Z[q^{\pm 1}]$. Let us denote by $S^{fr}$ the submodule of
$R{\cal K}^{fr}$ generated by framing expressions $K^{(1)}-qK$, for 
any framed knot in ${\cal K}^{fr}$. 
The quotient ${\cal S}_0(M,q)=R{\cal K}^{fr} /S^{fr}$ is called the 
framing skein module. 
\end{example}
We proved (with J.~Hoste) that if $M$ has no non-separating $S^2$  
then the module is free (with basis ${\cal K}$) \cite{H-P-2}
 and we were referred 
by D.~McCullough to his paper \cite{Mc} from which it follows that

\begin{theorem}\label{3.9}\ \\
 For a 3-manifold $M$, \ 
${\cal S}_0(M,q)= 
Z[q^{\pm 1}]{\cal K}^f\oplus \bigoplus_
{K \in {\cal K}-{\cal K}^f)} Z[q^{\pm 1}]/(q^{2}-1),$ 
where ${\cal K}^f$ is composed of knots which do not cut any 
$2$-sphere in $M$ in exactly one point.
\end{theorem}
Theorem 3.9 was also proved, by different method, by Chernov \cite{Cher}.

To summarize our examples from the geometric point of view, we stress 
that  a nonseparating $S^2$ in $M$ is detected by torsion in the 
skein module ${\cal S}_0^{fr}(M)$, a nonseparating $T^2$ in $M$ is 
detected by torsion in the
skein module ${\cal S}_{\pm}^{fr}(M)$, and a nonseparating oriented 
closed surface in $M$ is detected by torsion in the
skein module ${\cal S}_2^{fr}(M)$. More complicated skein modules 
detect, to various degree, separating surfaces as well (compare Section 8).

\section{Skein module based on Homflypt skein relation}\label{IX.4} 
\markboth{\hfil{\sc Skein modules}\hfil}
{\hfil{\sc Homflypt skein relation}\hfil}

The Jones-Conway (Homflypt) skein relation 
($v^{-1}L_+ - vL_- = z L_0$, where $L_+,L_-,L_0$ are oriented links
as in Fig. 4.1) is based on a relation which was hinted at by 
Alexander and Conway 
and can be thought of as a deformation of the crossing 
change (or, of a $2$-move; see Chapter VI). 

\begin{definition}\label{4.1} 
Let $M$ be an oriented 3-manifold, $R=Z[v^{\pm 1},z^{\pm 1}]$, $\cal L$ 
the set of all oriented links in $M$ up to ambient isotopy of $M$ and 
$S_3$ the submodule of $R\cal L$ generated by  the skein 
expressions $v^{-1}L_+ - vL_- -zL_0$. For convenience we allow the empty 
knot, $\emptyset$, and add the relation $v^{-1}\emptyset 
-v\emptyset -zT_1$, where $T_1$ denotes the trivial knot. 

Then the third skein module of $M$ 
is defined to be: $${\cal S}_3(M)={\cal S}_3(M;Z[v^{\pm 1},z^{\pm 1}], 
v^{-1}L_+-vL_--zL_0) = R{\cal L}/S_3 .$$ 
\end{definition} 
\ \\
\centerline{\psfig{figure=L+L-L0.eps}}
\centerline{Fig. 4.1}

The existence of the Homflypt 
polynomial of links in $S^3$ can be formulated as follows: 
\begin{theorem}[\cite{FYHLMO,P-T-1}] 
${\cal S}_3(S^3) = Z[v^{\pm 1},z^{\pm 1}]$, 
and the empty link is a generator of the module. 
In particular $L=P_L(v,z)T_1=P_L(v,z)(\frac{v^{-1}-v}{z})\emptyset$, 
where $P_L(v,z)$ is the Jones-Conway (Homflypt) polynomial of $L$. 
\end{theorem} 
The third skein module of a solid torus was computed by J.~Hoste and 
M.~Kidwell 
\cite{Ho-K} and independently by V.~Turaev \cite{Tu-2}. 
The first version of \cite{Ho-K}, from March 1987,
motivated me to define general skein modules a month later. 

The computations for $S^3$ and $S^1\times D^2$ are special case of the
general computation for the product of a surface and an interval, 
$F\times I$.  The third skein module of $F\times I$
has a structure of algebra ($L_1 \cdot L_2$ denotes the link obtained 
by placing $L_1$ above $L_2$) and the result is reminiscent of the 
Poincar\'e-Birkhoff-Witt theorem on the universal enveloping algebra of 
a Lie algebra.
\begin{theorem}[\cite{P-14}]
${\cal S}_3(F\times I)$ is an algebra which, as an $R$ 
module, is a free module isomorphic to the symmetric tensor algebra, 
${\bf S}R\hat\pi^o$, where $\hat\pi^o$ denotes the conjugacy classes of 
nontrivial elements of $\pi_1(F)$.
\end{theorem}

\begin{theorem}[\cite{Tu-3,P-15}]\label{IX.4.4} 
${\cal S}_3(F\times I)$ is an involutory Hopf algebra.
\end{theorem}

$S_3(F\times I)$  can be
interpreted as a quantization \cite{Ho-K,Tu-3,P-5,Tu-2,P-15}, and
${\cal S}_3(M)$ is related to the algebraic set of
$SL(n,\mathbb{C})$ representations of the
fundamental group of the manifold $M$, \cite{Si-4,Si-5}.

If $F$ is a non-orientable surface, then the twisted $I$-bundle 
over $F$ is an oriented 3-manifold ($M= F\hat\times I$) and we 
have, as in the product case, the concept of a diagram of a link.
This gives hope that the Homflypt skein module could also be computed 
in this case. M. Mroczkowski achieved this in the case 
of $F={\R}P^2$ proving, in particular, that the module is free \cite{Mro}. 

\begin{theorem}[\cite{Mro}]\label{IX.4.5}\ \\
${\cal S}_3(\mathbb{R}P^2 \hat\times I)$ is 
freely generated by standard oriented 
unlinks\footnote{The term unlink can be a little misleading as e.g. $lk(L_2)= 
\frac{1}{2}$ in $\mathbb{R}P^3$ 
(or $(\mathbb{R}P^2 \hat\times I)= \mathbb{R}P^3 \# D^3$), 
but we follow the terminology in \cite{Mro}.}
  $\vec L_n$. $\vec L_n$ is the link composed of $n$ copies 
of noncontractible simple closed curves on $RP^2$ as presented in 
Fig. 4.2 (${\R}P^2$ is presented as a disk with antipodal points identified).
$L_0$ is the empty link (if we allow it; otherwise, as in \cite{Mro}, 
one can take the trivial knot as the basic unlink $L_0$.).
\end{theorem}
\ \\
\centerline{\psfig{figure=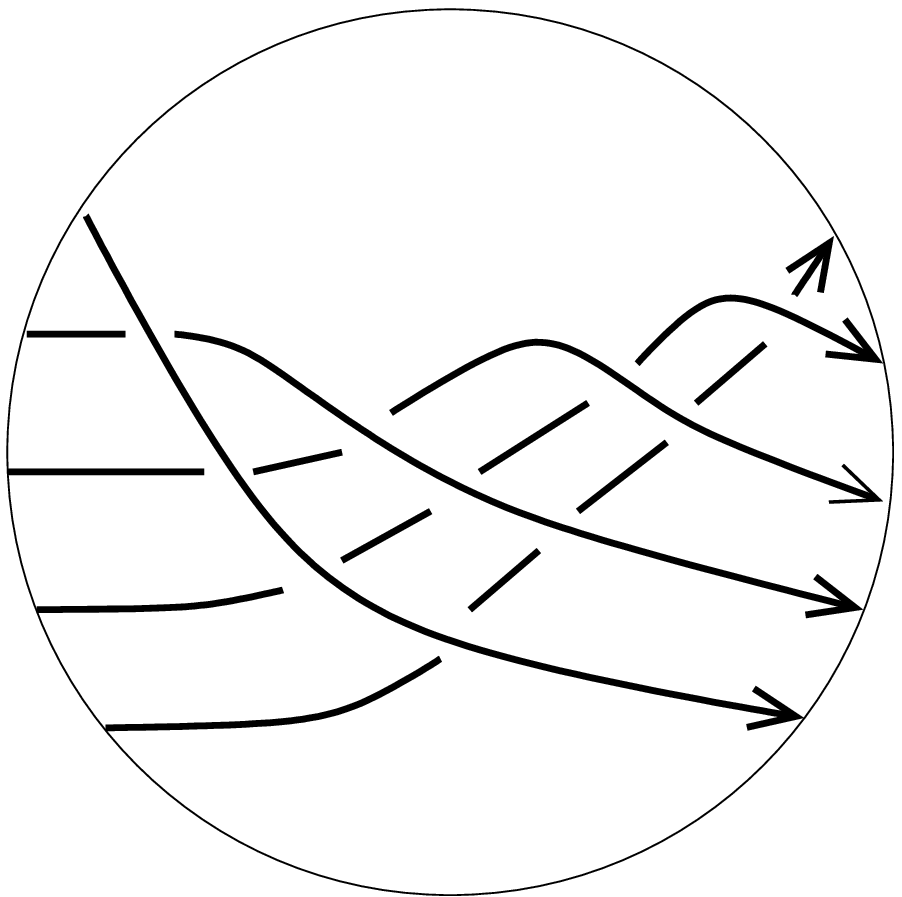,height=4.7cm}}
\begin{center}
Fig. 4.2; \ The standard oriented unlink  $\vec L_5$
\end{center}

Theorems 4.3 and 4.5 support the following conjecture
\begin{conjecture}\label{IX.4.6}\ \\ If $M$ is a submanifold of a 
rational homology sphere and it does not contain a closed, oriented 
incompressible surface then its Homflypt skein module ${\cal S}_3(M)$ 
is free and isomorphic to the symmetric tensor algebra over module 
spanned by conjugacy classes of nontrivial elements of the fundamental 
group, ${\cal S}_3(M) = {\bf S}R\hat\pi^o$.
\end{conjecture}

The idea of proving Conjecture 4.6 for lens spaces $L(n,1)$ 
is sketched in \cite{L-P}, but the paper is not finished yet.  

Skein modules can be defined for relative links (properly embedded 
1-dimensional manifolds) in the same manner as for absolute links 
(that is links composed only of closed curves). 
We assume usually that relative links are considered modulo ambient 
isotopy constant on the boundary. To define a relative skein module 
we should choose boundary points (or, at least, if the boundary 
is connected, the number of points\footnote{We give more detailed 
discussion of this point in the case of the relative Kauffman bracket 
skein module; see Section 7.}).  
As before, the accessible case is the product 
of a surface with boundary, $F$, and the interval. The case $S_3(F\times I,2n)$ 
for $F$ a disk is a skein module of tangles, closely related to an $A$-type 
Hecke algebra. It is known to be free of $n!$ generators.
 For $F$ being an annulus, we demonstrate after \cite{L-P}
 that the Homflypt skein module is free and related to a $B$-type 
Hecke algebra.
\begin{theorem}\label{4.7} Let $F=F_{0,2}$ be an annulus with $2n$ points, 
$p_1$,...,${p_{2n}}$ on one of its boundary components. We assume 
$p_1,...p_n$ are inputs and $p_{n+1},...,p_{2n}$ are outputs. 
${\cal S}_3(F\times I, 2n)$ denotes the relative Homflypt skein 
module with the chosen $2n$ points.
\begin{enumerate}
\item[(i)]
The module is an infinitely generated free $Z[v^{\pm 1},z^{\pm 1}]$-module.
\item[(ii)]
Let $R_0={\cal S}_3(F\times I)$ denote a ring which is the skein module of 
the product of an annulus and the interval. Then ${\cal S}_3(F\times I, 2n)$ 
is a free $R_0$ module with basis, $B(F_{0,2},2n)$
 described in subsequent definitions (\cite{L-P}).
\end{enumerate}
\end{theorem}

Sketch of the proof.\\
Our main tool is the Hoste-Kidwell-Turaev theorem (Theorem 4.3 for 
an annulus). First we need several definitions including the 
notion of ``half-way" Hecke algebra of type B, $H_n(p,q;\infty)$.
In what follows we think about $F_{0,2}\times I$ as the exterior of the 
(thickened) $z$-axis in $D^2 \times I$. The first, special string 
of a braid is represented by the $z$-axis
\begin{definition}\label{IX.4.8}\
\begin{enumerate}
\item[(i)]
The $B$-type Artin braid group $B_{1,n}= \{(t,
{\sigma}_1,...,{\sigma}_{n-1}\ \ | \ \ 
t{\sigma}_1t{\sigma}_1={\sigma}_1t{\sigma}_1t, t{\sigma}_i={\sigma}_it$
for $i> 1$, $ {\sigma}_i{\sigma}_j={\sigma}_j{\sigma}_i$ if $|j-i|>1$,
${\sigma}_i{\sigma}_{i+1}{\sigma}_i=
{\sigma}_{i+1}{\sigma}_i{\sigma}_{i+1}\}.$ It is a normal subgroup 
of a (type $A$) braid group $B_{n+1}$. The braid relation 
$t{\sigma}_1t{\sigma}_1={\sigma}_1t{\sigma}_1t$ is illustrated in 
Fig. 4.3. \\
\ \\
\centerline{\psfig{figure=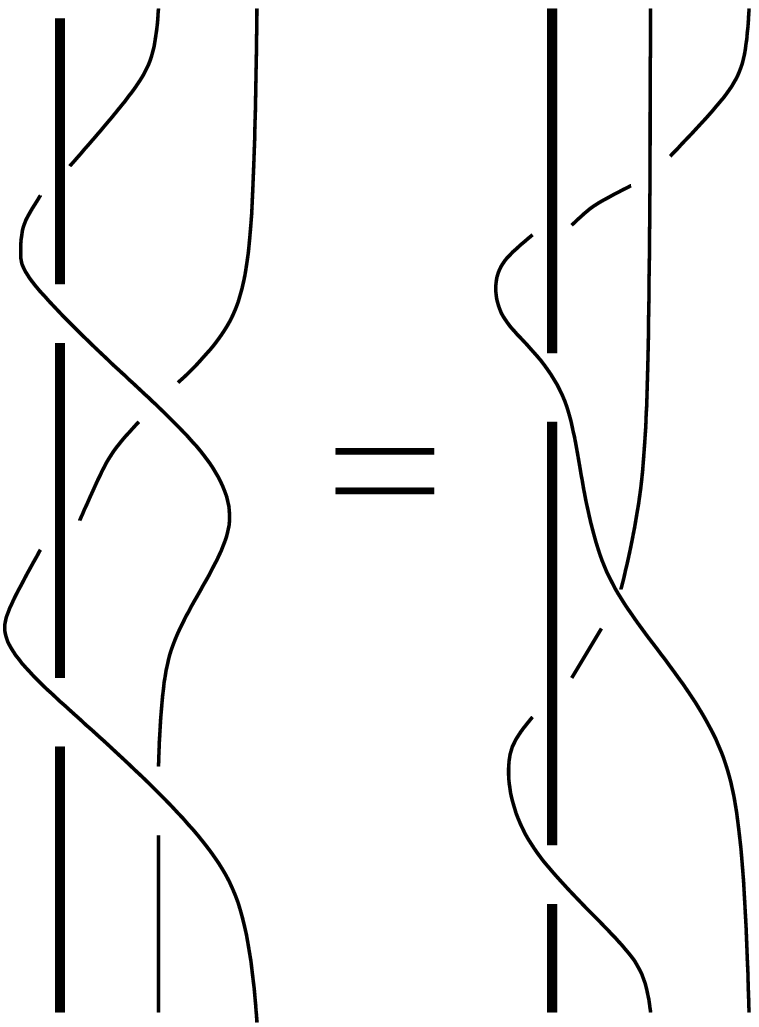,height=4.7cm}}
\centerline{Fig. 4.3}

\item[(ii)]
The ``half-way" Coxeter group
$W_n(\infty) = B_{1,n}/({\sigma}_i^2 = 1)= \{t,s_1,...,s_{n-1}:
ts_1ts_1=s_1ts_1t, ts_i=s_it$ if $i>1, s_is_j=s_js_i
$ if $ |i-j|>1$ and $s_is_{i+1}s_i=s_{i+1}s_is_{i+1} ,s_i^2 = 1\} $. 
We use $s_i$ for the image of $\sigma_i$ in the quotient space.
The group $W_n(\infty)$ (which is ``half-way" between Coxeter group 
and related Artin group because we do not assume the relation $t^2=1$) 
is well understood as it can be interpreted as
a group of weighted (or framed) permutations.
Thus we know that $W_n(\infty)  =
\mathbb{Z}^n \rhd\!\!\!\!< S_n$ (semidirect product) so we have
another familiar presentation coming from this
semi-direct product decomposition:
$W_n(\infty)=\{s_1,...s_{n-1},v_0,v_1,...,v_{n-1}\ | \ s_is_j=s_js_i
$ if $ |i-j|>1$, $s_is_{i+1}s_i= s_{i+1}s_is_{i+1},
s_i^2 = 1, v_iv_j=v_jv_i, s_iv_js_i = 
v_j$ for $j\neq i-1,j$, $s_iv_{i-1}s_i=v_i ,
s_iv_is_i=v_{i-1}\}$. 
When going from the second presentation to the first we have
$v_i=s_is_{i-1}...s_1ts_1...s_{i-1}s_i$ ($v_3$ is presented in Fig. 4.4). 
When we think of $W_n(\infty)$ as of ``framed"
permutations then we interpret $v_i$ as one positive twist on the framing of
the $i+1$ string. 
\\
\ \\
\centerline{\psfig{figure=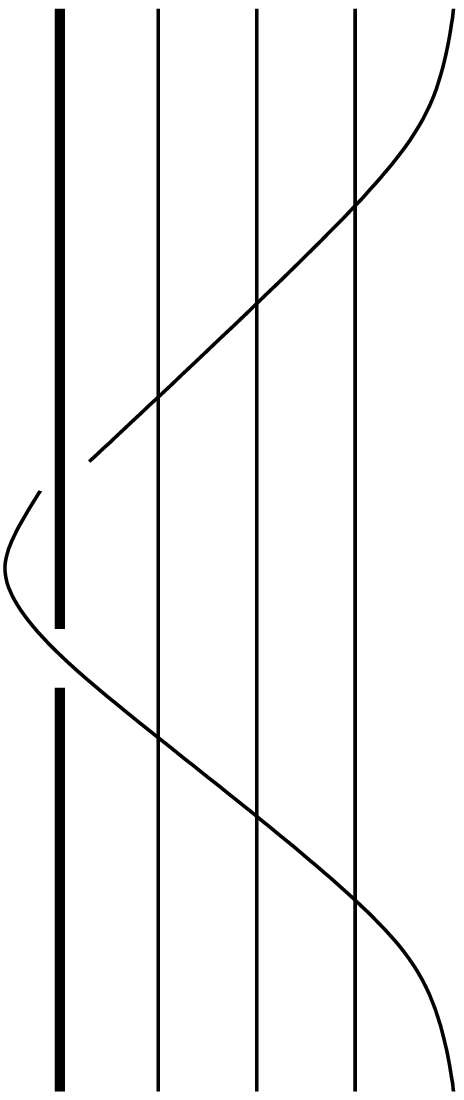,height=4.7cm}}
\centerline{Fig. 4.4}

\item[(iii)] For $R=Z[p^{\pm 1},q^{\pm 1}]$ we define 
$H(p,q;\infty)$ (the ``half-way" Hecke algebra) 
as the quotient of the group ring $RB_{1,n}$
by the quadratic relations: ${\sigma}_i^2=p{\sigma}_i+q$. 
We can write, therefore,
the presentation of $H(p,q;\infty)$ as follows:\\
$H_n(p,q;\infty) = \{t,g_1,...,g_{n-1}\ | \
 tg_1tg_1=g_1tg_1t,\ tg_i=g_it\ if i>1,\
 g_ig_j=g_jg_i$
if $|i-j|>1 $ and $g_ig_{i+1}g_i=g_{i+1}g_ig_{i+1} ,g_i^2 = pg_i + q\} .$\\
We use $g_i$ for the class of $\sigma_i$.
The ``half-way" Hecke algebra $H(p,q;\infty)$ 
is a deformation of $RW_n(\infty)$
with respect to the equations $s_i^2=1$, where $R=Z[p^{\pm 1},q^{\pm 1}]$ 
and $RW_n(\infty)$ is the group algebra of $W_n(\infty)$.
\end{enumerate}
\end{definition} 
\begin{proposition}[Natural generating set; basis]\label{IX.4.9}\
\begin{enumerate}
\item[(i)] Every element of $W_n(\infty)$ can be uniquely written in the 
normal form. We give an inductive definition. 
For $n=1$ the normal form is $t^i$, $i\in \mathbb{Z}$. 
If normal words for $W_n(\infty)$ are defined, 
then normal words for $W_{n+1}(\infty)$ are of the form
$s_ks_{k-1}...s_1t^js_1...s_nw$, $j\neq 0$ or $s_{i}s_{i+1}...s_{n}w$ 
where $w$ is a normal word in $W_n(\infty)$ and $k\leq n$, $1\leq i \leq n+1$. 
 In particular these words 
form a basis of $RW_n(\infty)$.
\item[(ii)] The words from (i) where we take $g_i$ in place of $s_1$ 
form a generating set of $H(p,q;\infty)$. In Figure 4.5 we illustrate 
$g_2g_1t^2g_1g_2g_3g_4g_1t^{-1}g_1g_2g_3g_1g_2g_1$. 
\item[(iii)] The generating set described in (ii) is 
a basis of $H(p,q;\infty)$.
\end{enumerate}
\end{proposition}
\ \\
\centerline{\psfig{figure=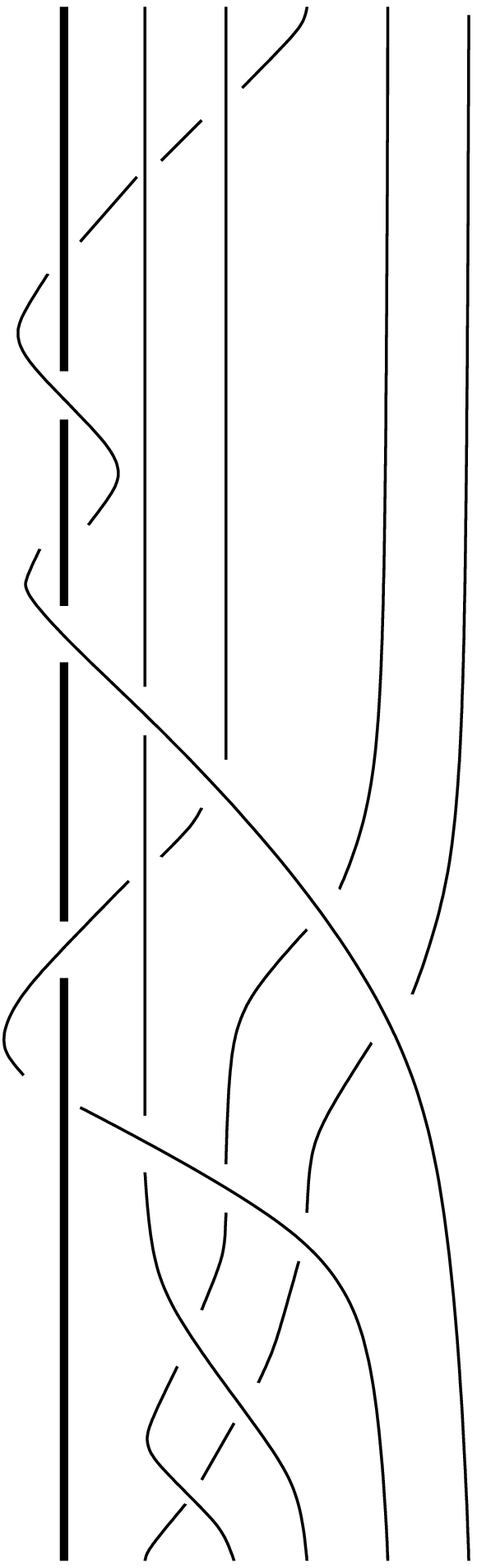,height=4.7cm}}
\centerline{Fig. 4.5}
Sketch of the proof.

Part (i) follows easily from our ``framed" permutation interpretation 
of $RW_n(\infty)$. Part (ii) can be proved by induction on the number 
of letters in a braid (quadratic relations allows us to reduce the number of 
letters (crossings) or move us closer to a normal form). To prove that 
our generating set is in fact a basis of $H(p,q;\infty)$ is more 
difficult but we can do it very quickly if we assume the theorem of 
Kidwell-Hoste and Turaev about the structure of the skein
module of the solid torus. It is the same trick used
by Morton and Traczyk \cite{Mo-Tr-2} to find the skein module of a
tangle when the existence of the Homflypt polynomial of links in $S^3$ 
is established.

We show that our generators from (ii) are linearly independent by 
showing that certain form is nondegenerated.
It is more convenient here to consider a slightly adjusted version 
of the ``half-way" Hecke algebra, $H_n(v,z;\infty)$ 
with the relation $v^{-1}s_i^2=zs_i+v$, which is also a Homflypt skein 
module relation.
Let $e_1,e_2,...$ denote our generating set of $H_n(v,z;\infty)$ 
from (ii).  We prove that the above elements are linearly independent.
Let $\sum_i a_i(v,z)e_i=0$ and $j$ be an index in the sum for which $z$ 
degree is the smallest.
Consider the bilinear function:
$$\Phi:H_n(v,z;\infty) \times H_n(v,z;\infty) \to {\cal S}_3(S^1\times D^2)$$
Where $\Phi (\alpha,\beta)$ is obtained by composing $\alpha$ and $\beta$
first, then closing the result and filling up the braid 
axis\footnote{Equivalently we consider braids up to conjugation 
and Markov moves.}. The result
is a link in the solid torus (with given $I$-bundle structure) and we
consider its value in ${\cal S}_3(S^1\times D^2)$ as $\Phi (\alpha,\beta)$.
Now consider $\Phi (\sum_i a_i(v,z)e_i, e_j^{-1})$, in the standard
basis of $ {\cal S}_3(S^1\times D^2)$, composed of layered families of 
torus knots of type ($2,n$). We can see easily that the unique
coefficient of the basis element with the lowest power of $z$ is given by
$\Phi (a_j(v,z)e_j,e_j^{-1})= a_j(v,z)(v^{-1}-v)^nz^{-n}$. Thus from
$\sum_i a_i(v,z)e_i=0$ follows $a_j(v,z)=0$ and $a_i(v,z)=0$ for any $i$.
These complete the proof of Proposition 4.9.

To complete a proof of Theorem 4.7 notice 
that in ${\cal S}_3(F_{0,2}\times I,2n)$ we 
start from $n$ arcs and some number of closed components in 
$F_{0,2}\times I$ and we assume that arcs have orientations  
from the level $F_{0,2}\times \{1\}$ to $F_{0,2}\times \{0\}$ 
(we call such relative links -- string links).
Using Homflypt skein relations we can reduce every string link 
to a linear combination of elements of the
generating set from Proposition 4.9(ii)  (over the ring $R_0$ 
as we have to allow closed curves). We organize our reduction 
along some properly chosen complexity: for arcs, distance from 
the normal form, and for closed curves, how far they are from the 
basic curves of ${\cal S}_3(F_{0,2}\times I)$.
The fact that our braids from (ii) are linearly independent allows us to 
complete the proof of Theorem 4.7.

As mentioned at the end of Section 3, torsion in skein modules is 
often related to surfaces in a 3-manifold. We illustrate this for the 
Homflypt polynomial in one simple instance: a nonseparating 2-sphere.
\begin{proposition}\label{IX.4.10}
Let $S^2$ be a nonseparating 2-sphere in an oriented 3-manifold $M$ and 
$L$ a link cutting $S^2$ in exactly one point. Then 
$((\frac{v^{-1}-v}{z})^2-1)L=0$ in ${\cal S}_3(M)$.
\end{proposition} 
\begin{proof} 
Consider the link $L'$ composed of $L$ and meridian circle, $C$, around it 
(Hopf component), Fig. 4.6. On one hand $C$ can be taken out of $L$ 
using the ``other side" of $S^2$ (thus $C$ is a trivial circle in $M-L$).
Therefore $L'= \frac{v^{-1}-v}{z}L$ in ${\cal S}_3(M)$. 
On the other hand we can use skein relation to compute $L'$ and 
get $L'= v^2(\frac{v^{-1}-v}{z})L + vzL$, Fig. 4.6. Therefore 
$\frac{v^{-1}-v}{z}L = (v^2(\frac{v^{-1}-v}{z}) + vz)L$ and 
$(frac{v^{-1}-v}{z}(v^2-1) + vz)L=0$, and finally 
$((\frac{v^{-1}-v}{z})^2-1)L=0$.
\end{proof}

\ \\
\centerline{\psfig{figure=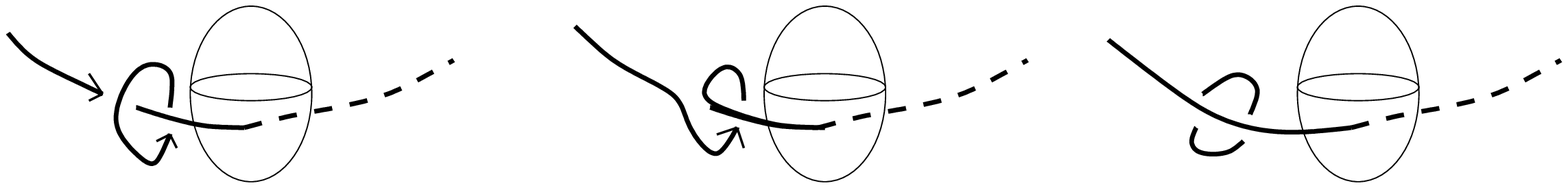,height=1.6cm}}
\centerline{Fig. 4.6}

\section{Homotopy skein module}\label{5} 
If we ignore self-crossings in the third (Homflypt) 
skein module we obtain its simplified 
version which we call the homotopy skein module \cite{H-P-2}. 

\begin{definition}\label{5.1} 
Let $M$ be an oriented 3-manifold, $\cal L$ the set of oriented links in $M$,
$R$ a commutative ring with unit, and $z$ a fixed element of $R$. The homotopy
skein module ${\cal H}{\cal S}(M,R,z)$ is defined as the quotient of the free
$R$ module over $\cal L$ and the submodule ${\cal H}$ generated by skein
relations $L_+-L_-$ for a selfcrossing and $L_+- L_- - zL_0$ for a crossing
between different components of $L_{\pm}$. That is,
${\cal H}{\cal S}(M,R,z)= R{\cal L}/{\cal H}$.
\end{definition}

If $M=F \times [0,1]$, where $F$ is an oriented surface, then we have a
multiplication in the module: $L_1 \cdot L_2$ denotes the link obtained
by placed $L_1$ above $L_2$ in $F \times [0,1]$ with respect to ``height" 
$[0,1]$.  If we assume that $\cal L$ contains the empty knot then 
${\cal H}{\cal S}(M;R,z)$ becomes an algebra (in \cite{H-P-2} a skein 
module with the empty set allowed is called a reduced skein module).
\begin{theorem}[\cite{H-P-2,Tu-2}]\ 
\begin{enumerate}
\item[(a)] ${\cal H}{\cal S}(F\times [0,1];R,z)$ is an algebra which as an $R$
module is a free module isomorphic to the symmetric tensor algebra,
${\bf S}R\hat\pi$, where $\hat\pi$ denotes the conjugacy classes of elements
of $\pi_1(F)$.
\item[(b)] ${\cal H}{\cal S}(F\times [0,1];R,1)$ is an $R$-algebra isomorphic
to the universal enveloping algebra of the Goldman-Wolpert Lie
algebra $U(R\hat\pi)$. Going ``backwards" we can describe the Lie bracket
$[\alpha,\beta]$ on $R\hat\pi$ as a projection on $F$ of 
$K_{\alpha}\cdot K_{\beta} - K_{\beta} \cdot K_{\alpha}$, 
where $K_{\alpha}$ and $K_{\beta}$ are knots representing curves 
$\alpha$ and $\beta$ respectively.
\item[(c)] ${\cal H}{\cal S}(F\times [0,1];R,z)$ is a Drinfeld-Turaev
quantization of the Goldman-Wolpert Poisson algebra of curves on $F$.
\end{enumerate}
\end{theorem}
The following $q$-analogue of the homotopy skein module is of considerable 
interest. 

\begin{definition}
Let $M$ be an oriented 3-manifold, $\cal L$ the set of oriented links in $M$,
$R$ a commutative ring with unit, $z$ a fixed element of $R$, and $q$
an invertible element of $R$. The $q$ homotopy
skein module ${\cal H}^q{\cal S}(M;R,q,z)$ is defined as the quotient of the 
free $R$ module over $\cal L$ and the submodule ${\cal H}^q$ generated by skein
relations $L_+-L_-$ for a selfcrossing and $q^{-1}L_+-qL_- - zL_0$ for a 
crossing between different components of $L_{\pm}$. That is,
${\cal H}{\cal S}^q(M;R,q,z)= R{\cal L}/{\cal H}^q$.
\end{definition}
We use the notation ${\cal H}{\cal S}^q(M)$ for 
${\cal H}{\cal S}^q(M;Z[q^{\pm 1},z],q,z)$. For $M=S^3$ the skein module
is freely generated by trivial links and the value of a link in the
module depends only on the number of components and linking numbers between 
components (linking matrix) \cite{P-6,P-12}. The skein module 
can be described using the dichromatic polynomial of the associated 
weighted graph:
vertices correspond to link components and weights of edges correspond to
linking numbers. An interesting feature of the skein module is that it
distinguishes some links with the same Jones-Conway (Homflypt)
polynomial.

\begin{example}\label{IX.5.4}\ \\
The 3-component links in Fig. 5.1 have the same Jones-Conway 
polynomial \cite{Bi-2} but they are different in ${\cal H}{\cal S}^q(M)$.
Namely $L_1= -(1+q^2+q^4+q^6)z^2T_1 + (q^{-1}+q+q^3+q^5-q^7)zT_2 + q^{6}T_3$
and $L_2= -(q^4+2q^6+q^8)z^2T_1 + (q+2q^3+2q^5-q^7-q^9)zT_2 + q^{6}T_3$. 
\end{example}
\ \\
\centerline{\psfig{figure=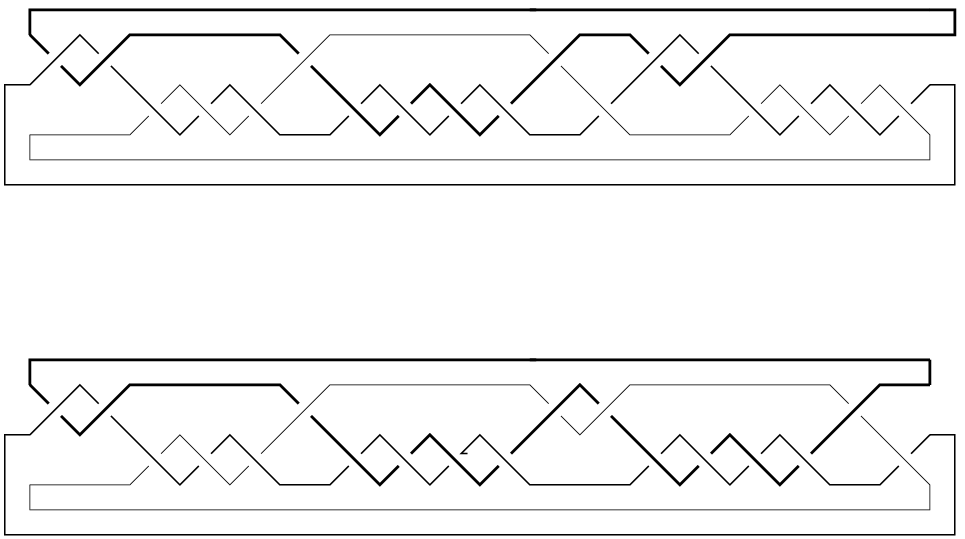}} 
\begin{center} 
Fig. 5.1. 
\end{center} 
We generalize, partially the case of $M=S^3$ or rather $M=D^3=D^2\times I$ 
into the case of $M=F\times I$.
Let ${\cal L}^{h}$ denote the set of homotopy links in $M$, that is,
${\cal L}^{h} = {\cal L}/(L_{+} - L_{-})$ where relations are
yielded by self-crossings. Let $\hat \pi$ denote
the set of conjugacy classes in $\pi_{1} (M)$, or equivalently the
set of homotopy knots in $M$. Choose some linear ordering, denoted
by $\leq$, of elements of $\hat \pi$. Given a homotopy link
$L=\{K_{1},K_{2},\ldots ,K_{n}\}$ in $F\times I$, we shall say
that $L$ is a layered homotopy link with respect to the ordering
of $\hat \pi$ if each $K_{i}$ is above $K_{i+1}$ in $F\times I$
and $K_{i}\leq K_{i+1}$. Let ${\cal B}$ be the set of all layered
homotopy links with respect to the ordering of $\hat \pi$,
including the empty link.
\begin{theorem}\label{IX.5.5}\ 
\begin{enumerate}
\item[(i)] The $q$-homotopy skein module ${\cal H \cal S}^{q}(F\times I)$
is generated by ${\cal B}$.
\item[(ii)] The homotopy skein module ${\cal H \cal S}(F\times I)$
is freely generated by ${\cal B}$; \cite{H-P-1}.
\item[(iii)] If $\pi_{1} (F)$ is abelian then the $q$-homotopy skein
module ${\cal H \cal S}^{q}(F\times I)$ is freely generated by
${\cal B}$.
\end{enumerate}
\end{theorem}
For the proof we refer to \cite{P-6,P-30}.

If the Euler characteristic, $\chi (F)$, is negative then the
$q$-homotopy skein module has torsion.
 This is described in Theorem 5.6. On the other hand we
know,  Theorem 5.2(a), that for $q=\pm 1$ the module is free.

\begin{theorem}\label{IX.5.6}
Let $F$ be a surface (not necessary compact) which contains a disc
with 2 holes or a torus with a hole embedded $\pi_1$-injectively;
equivalently, $\pi_1 (F_{0})$ is not abelian for a connected
component $F_{0}$ of $F$ (in the compact connected case this
 means that $\chi (F)<0$). Then
\begin{enumerate}
\item[(a)]
${\cal H \cal S}^q(F\times I)$ has torsion.
\item[(b)] Let $\alpha : R{\cal B} \to {\cal H \cal S}^q(F\times I)$
be an $R$-homomorphism given by $\alpha (L)=L$. Then
$ker \alpha \neq \{0\}$.
\end{enumerate}
\end{theorem}

U. Kaiser generalized Theorem 5.6, fully characterizing oriented 
3-manifolds for which $q$-homotopy skein module has torsion \cite{Kai-2}.

\section{The Kauffman bracket skein module}\label{IX.6}

The skein module based on the Kauffman bracket relation is, so far,
the most extensively studied object in {\it algebraic topology based 
on knots}.
We describe in this section the basic properties of the Kauffman Bracket
Skein Module (KBSM) and list manifolds for which the structure of 
the module is known. In the seventh section, we give the detailed proof
of the structure of KBSM of a 3-manifold which is an interval bundle over
a surface. In the ninth section we analyze the structure of algebra
 in the case of a surface times an interval, we introduce the notion
of the skein algebra of a group and investigate the relation with
representations of the group to $SL(2,C)$.

\begin{definition}[\cite{P-5,H-P-3}]\label{6.1}\ \\
Let $M$ be an oriented 3-manifold, ${\cal L}_{fr}$ the set of unoriented
framed links in $M$ (including the empty knot, $\emptyset$), 
$R$ any commutative ring with identity and $A$ its
invertible element. Let $S_{2,\infty}$ be
the submodule of $R{\cal L}_{fr}$ generated by skein expressions
$L_+- AL_0 - A^{-1}L_{\infty}$, where the triple 
$L+, L_0, L_{\infty}$ is presented in Fig.6.1, and $L \sqcup T_1 +(A^2 +
A^{-2})L$, where $T_1$ denotes the trivial framed knot. 
We define the Kauffman bracket skein module, ${\cal S}_{2,\infty}(M;R,A)$,  
as the quotient ${\cal S}_{2,\infty}(M;R,A)= R{\cal L}_{fr}/S_{2,\infty}$.
\end{definition}
\ \\
\centerline{\psfig{figure=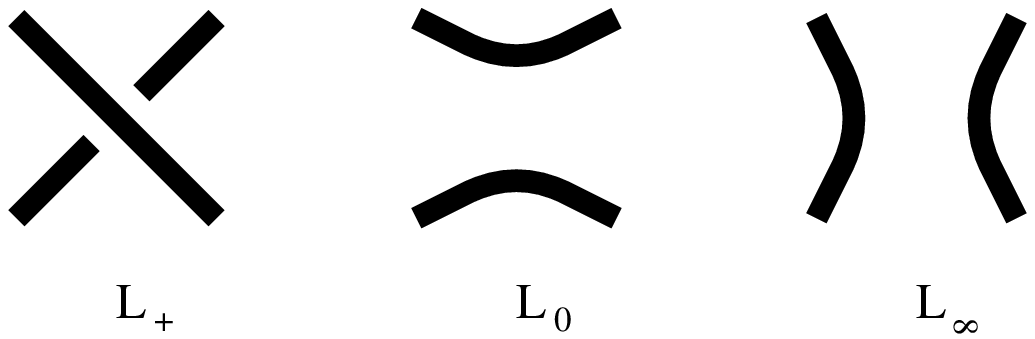}}
\begin{center}
Fig. 6.1.
\end{center}
Notice that $L^{(1)}=-A^3L$ in ${\cal S}_{2,\infty}(M;R,A)$; we call this the
framing relation. 
We use the simplified notation 
${\cal S}_{2,\infty}(M)$ for ${\cal S}_{2,\infty}(M;Z[A^{\pm 1}],A)$.

We list below several elementary properties of KBSM including
description of the 
KBSM of any compact 3-manifold using generators and relators. 
\begin{proposition} \
\begin{enumerate} 
\item [(1)] An orientation preserving embedding of 3-manifolds 
$i: M \to N$ yields the homomorphism of skein modules 
$i_*: {\cal S}_{2,\infty}(M;R,A) \to {\cal S}_{2,\infty}(N;R,A)$. 
The above correspondence leads to a functor from the category of 
3-manifolds and orientation preserving embeddings (up to ambient isotopy) 
to the category of $R$-modules (with a specified invertible element $A \in R$). 
\item [(2)] 
\begin{enumerate} 
\item [(i)] 
If $N$ is obtained from $M$ by adding a 3-handle to it (i.e. capping off 
a hole), and $i: M \to N$ is the associated embedding, 
then $i_*:{\cal S}_{2,\infty}(M;R,A) \to {\cal S}_{2,\infty}(N;R,A)$ 
is an isomorphism. 
\item [(ii)] If $N$ is obtained from $M$ by adding a 2-handle to it, 
and $i: M \to N$ is the associated embedding, 
then $i_*:{\cal S}_{2,\infty}(M;R,A) \to {\cal S}_{2,\infty}(N;R,A)$ 
is an epimorphism. 
\end{enumerate} 
\item [(3)] If $M_1 \sqcup M_2$ is the disjoint sum of 3-manifolds $M_1$ 
and $M_2$ then 
$${\cal S}_{2,\infty}(M_1\sqcup M_2;R,A)= {\cal S}_{2,\infty}(M_1;R,A) 
\otimes {\cal S}_{2,\infty}(M_2;R,A).$$ 
\item [(4)] 
(Universal Coefficient Property)\\ 
Let $r: R \to R'$ be a homomorphism of rings (commutative with 1). 
We can think of $R'$ as an $R$ module. Then the identity map on 
${\cal L}_{fr}$ induces the isomorphism of $R'$ (and $R$) modules: 
$$ \bar r: {\cal S}_{2,\infty}(M;R,A)\otimes_R R' \to 
{\cal S}_{2,\infty}(M;R',r(A))  .$$ 
\item [(5)] 
Let $(M,\partial M)$ be a 3-manifold with the boundary $\partial M$, 
and let $\gamma $ be a simple closed curve on the boundary. Let $N=M_{\gamma}$ 
be the 3-manifold obtained from $M$ by adding 
a 2-handle along $\gamma $. Furthermore let ${\cal L}_{fr}^{gen}$ be a set 
of framed links in $M$ generating ${\cal S}_{2,\infty}(M;R,A)$.\\ 
Then ${\cal S}_{2,\infty}(N;R,A) = {\cal S}_{2,\infty}(M;R,A)/J$, where 
$J$ is the submodule of ${\cal S}_{2,\infty}(M;R,A)$ generated by 
expressions $L- sl_{\gamma}(L)$, where $L\in {\cal L}_{fr}^{gen}$ and 
$sl_{\gamma}(L)$ is obtained from $L$ by sliding it along ${\gamma}$ 
(i.e. handle sliding). 
\item [(6)] 
Let $M$ be an oriented compact manifold and consider its Heegaard 
decomposition (that is $M$ is obtained from the handlebody $H_n$ 
by adding 2- and 3-handles to it, then $M$ has a presentation 
as follows: generators of 
${\cal S}_{2,\infty}(M;R,A)$ are generators of ${\cal S}_{2,\infty}(H_n;R,A)$ 
and relators are yielded by 2-handle slidings. 
\end{enumerate} 
\end{proposition} 
\begin{proof} 
\begin{enumerate} 
\item [(1)] $i_*$ is well defined because if framed links $L_1$ and $L_2$
are ambient isotopic in $M$ then $i(L_1)$ and $i(L_2)$ are ambient 
isotopic in $N$. Furthermore any skein triple $L_+,L_0,L_{\infty}$ in $M$, 
is sent by $i$ to a skein triple in $N$. Finally $i(T_1)$ is a trivial
framed knot in $N$.\ Notice that if $i_*:M \to N$ is an orientation reversing 
embedding then $i_*$ is a $Z$-homomorphism and $i(Aw)=A^{-1}i(w)$. 
\item [(2)] 
\begin{enumerate} 
\item [(i)] It holds because the cocore of a 3-handle 
is $0$-dimensional.\footnote{A manifold $N$ is obtained 
from an $n$-dimensional manifold $M$
by attaching a $p$-handle, $D^p \times D^{n-p}$, to $M$, if $N= M \cup_f
D^p \times D^{n-p}$ where $f: \partial D^p \times D^{n-p}$ is an embedding.
$D^p \times \{0\}$ is a core of the handle and 
$\{0\} \times D^{n-p}$ is a cocore of the handle \cite{R-S}.} 
\item [(ii)] It holds because the cocore of a 2-handle is $1$-dimensional. 
\end{enumerate} 
\item [(3)] This is a consequence of the well known property of short exact 
sequences, \cite{Bl}:\\ 
If $ 0 \to A' \to A \to A'' \to 0$ and $ 0 \to B' \to B \to B'' \to 0$ are 
short exact sequences of $R$-modules then 
$ 0 \to A'\otimes B + A\otimes B' \to A\otimes B \to A''\otimes B" \to 0$ 
is a short exact sequence. 
\item [(4)] 
The exact sequence of $R$ modules 
$$ S_{2,\infty}(R,A) \to R{\cal L}_{fr} \to {\cal S}_{2,\infty}(M;R,A) \to 0$$ 
leads to the exact sequence of $R'$ modules (\cite {C-E}, Proposition 4.5): 
$$ S_{2,\infty}(R,A)\otimes_R R' \to R{\cal L}_{fr}\otimes_R R' 
\to {\cal S}_{2,\infty}(M;R,A) \otimes_R R' \to 0.$$ 
Now, applying the "five lemma" to the commutative diagram with exact rows 
(see for example \cite{C-E} Proposition 1.1) 
\begin{displaymath} 
\begin{array}{ccccccc} 
 S_{2,\infty}(R,A) \otimes_R R' & \to & R{\cal L}_{fr}\otimes_R R' 
& \to & {\cal S}_{2,\infty}(M;R,A) \otimes_R R' & \to & 0 \\ 
\downarrow epi & & \downarrow iso & & \downarrow \bar r & & \\ 
S_{2,\infty}(R',\bar r (A)) & \to & R'{\cal L}_{fr} & \to & 
{\cal S}_{2,\infty}(M;R',\bar r(A))  & \to & 0 
\end{array} 
\end{displaymath} 
we conclude that $\bar r$ is an isomorphism of $R'$ (and $R$) modules. 
\item[(5)]
It follows from (2)(ii) because any skein relation can be performed in $M$,
and the only difference between KBSM of $M$ and $N$ lies in the fact that 
some nonequivalent links in $M$ can be equivalent in $N$; 
the difference lies exactly in the possibility 
of sliding a link in $M$ along the added 2-handle (that is $L$ is moving 
from one side of the cocore of the 2-handle to another).
\item[(6)] It follows from (5) and (2)(i).
\end{enumerate} 
\end{proof} 

\begin{remark}\label{IX.6.3}
The Universal Coefficient Property holds for all skein modules considered
in this work. It is the case because we always have the isomorphism
$R{\cal L}{\otimes}_R R' \to R'{\cal L}$ and the epimorphism of
submodules of skein relations $S(R){\otimes}_R R' \to S(R')$, thus
the proof described in (5), based on ``five lemma," works for our
skein modules.
\end{remark}
In the next theorem we list manifolds for which the exact structure of
the Kauffman bracket skein module has been computed.
\begin{theorem}[\cite{P-5,H-P-4,H-P-5,H-P-6,Bu-2}]\label {6.4}\ \\
\begin{enumerate}
\item
[(a)]  ${\cal S}_{2,\infty}(S^3)=Z[A^{\pm 1}]$, 
more precisely: $\emptyset$ is the generator of the module
and $L=<L> T_1=(-A^2-A^{-2}<L> \emptyset$ where $<L>$ is the Kauffman 
bracket polynomial of a framed link $L$.
\item
[(b)]  ${\cal S}_{2,\infty}(F\times [0,1])$ is a free module generated
by links (simple closed curves) on  $F$ with no trivial component 
(but including the empty knot). 
Here $F$ denotes  an oriented surface (see also Theorems 5.1 and 5.9).\\ 
This applies in particular to a handlebody, because 
$H_n = P_n \times I$, where $H_n$ is a handlebody of genus $n$ and $P_n$ is 
a disc with $n$ holes. 
\item
[(c)]  ${\cal S}_{2,\infty}(L(p,q))$ is a free $R$ module and it has
$[p/2]+1$ generators, where $[x]$ denotes the integer part of $x$. 
\item
[(d)] ${\cal S}_{2,\infty}(S^1\times S^2)=  Z[A^{\pm 1}]\oplus 
\bigoplus_{i=1}^{\infty} Z[A^{\pm 1}]/(1-A^{2i+4})$
\item
[(e)] The skein module of the complement of the torus knot
of type $(k,2)$ is free.
\item 
[(f)] Let $W$ be the classical Whitehead manifold, then ${\cal S}_{2,\infty}(W)$
is infinitely generated torsion free but not free.
\end{enumerate} 
\end{theorem}
In \cite{H-P-5}, (f) is proved for a large class of genus one Whitehead
type manifolds, and the paper ends with an optimistic note: {\it for 
the classical Whitehead manifold it seems
feasible to find the exact structure of ${\cal S}_{2,\infty}(W)$, 
and we plan to address that in a future paper.}

We prove (b), with its generalizations, in the next section. 
The part (e) was generalized recently by Thang Le (\cite{Le} April 2004) 
to  the Kauffman bracket skein modules of the
 exteriors of 2-bridge knots.
\begin{theorem}[T.Le]\label{6.5}
For a 2-bridge (rational) knot $K_{\frac{p}{m}}$ the
skein module is the free $\mathbb{Z}[A^{\pm 1}]$ module with the basis
$\{x^iy^j\}$, $0\leq i$, $0\leq j \leq \frac{p-1}{2}$, where $x^iy^j$
denotes the element of the skein module represented by the link
composed of $i$ parallel copies of the meridian curve $x$ and $j$
parallel copies of the curve $y$; see Fig. 6.2.
Le's theorem generalizes results in \cite{Bu-3} and \cite{B-L}.
\end{theorem}
\ \\
\ \\
\centerline{\psfig{figure=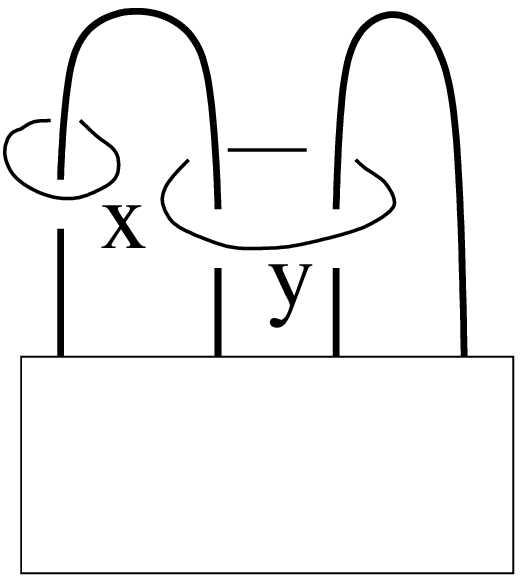,height=5.7cm}}
\centerline{Fig. 6.2}

One should compare here Le's Theorem with the part (c) of   
computation of Hoste-Przytycki \cite{H-P-4}. These two results are 
related by the fact that the lens space $L(p,m)$ is the double 
branched cover of $S^3$ branched along the 2-bridge 
knot $K_{\frac{p}{m}}$. We will discuss this in 
the monograph \cite{P-30} . 
 \\
\ \\
\newpage

\section{KBSM and relative KBSM of $F\times I$ and 
$F\hat{\times} I$ }\label{IX.7}\ \\
\markboth{\hfil{\sc Skein modules}\hfil}
{\hfil{\sc KBSM and relative KBSM of $F\times I$}\hfil}

The understanding of the Kauffman bracket skein module of the product
of a surface and the interval is the first step to understanding KBSM
of a general 3-manifold. Furthermore the case of $F\times I$ is 
relatively easy to understand because we can project links onto the surface
and work with diagrams of links. This can be generalized to twisted 
$I$-bundles over $F$ and one can have reasonable hopes that the method
can work for other 3-manifolds by projecting links to spines of 3-manifolds.
The relative case is described in Theorem 7.10.

\begin{theorem}\label{7.1}
Let $M$ be an oriented 3-manifold which is either equal to $F\times I$,
where $F$ is an oriented surface, or it equal to a twisted $I$ bundle over
$F$ ($F\hat{\times} I$), where $F$ is an unorientable surface. Then the KBSM, 
${\cal S}_{2,\infty}(M;R,A)$, is a free $R$-module with 
a basis $B(F)$ consisting of links in $F$ without contractible 
components (but including the empty knot). 
\end{theorem}

\begin{proof} We will give here the proof of Theorem 5.1 which is based on the
original proof of Kauffman on the existence of his bracket polynomial.\\
Let $M$ be an oriented 3-manifold which is an $I$-bundle over a 
surface $F$.\footnote{Because $M$ is 
oriented therefore for $\gamma$ in $F$ changing 
orientation of $F$, the restriction of the $I$-bundle to $\gamma$ is a
nontrivial bundle (M\"obius band). For $\gamma$ preserving orientation 
of $F$, the bundle is trivial (an annulus).}
Let $B(F)$ consist of all links in $F$ which have no trivial components
(including $\emptyset$). Furthermore each link is equipped with
an arbitrary, but specific framing (to be concrete we can assume that if
a knot in $F$ preserves the orientation of $F$ then we choose as its framing
the regular neighborhood of $K$ in $F$ (``blackboard" framing), if $K$ is
changing the orientation on $F$ then its regular neighborhood is a M\"obius
band so to get a framing we perform a positive half twist on it). 
Now one can quickly see that $B(F)$ is a generating set of 
${\cal S}_{2,\infty}(M;R,A)$. 
Namely every link
in $M$ has a regular projection on $F$ and any link can be reduced by  skein 
relations so that a projection has no crossings. 
Then another relation allows us to eliminate trivial components and finally 
the framing relations allow us to adjust framing. 
We will prove that $B(F)$ is a basis for 
${\cal S}_{2,\infty}(M;R,A)$. First we need however to consider the space of 
link diagrams (for a nonorientable surface $F$ the proof is still easy but 
requires great care).
\newpage

\begin{definition}\label{7.2}\ 
\begin{enumerate}
\item[(a)]
A link diagram on $F$ is a 4-valent graph in $F$ 
(allowing loops without vertices) such that one corner of each vertex 
is marked. $F$ does not need to be oriented for this definition. 
\item[(b)] 
Let $\cal D$ be a set of link diagrams on $F$ (up to isotopy of $F$),
and $R\cal D$ the free module over $\cal D$. The skein space of
diagrams, ${\cal S}\cal D$ is defined as a quotient:
$${\cal S}{\cal D}(F;R,A)= R{\cal D}/(
\parbox{0.5cm}{\psfig{figure=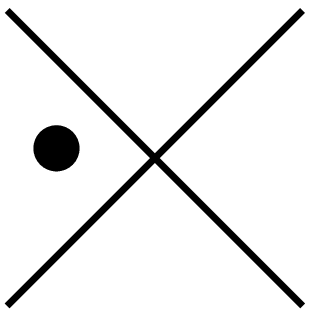,height=0.4cm}} -
A\parbox{0.5cm}{\psfig{figure=L0nmaly.eps}} - 
A^{-1} \parbox{0.5cm}{\psfig{figure=Linftynmaly.eps}},\ \
 D\sqcup T_1 +(A^2 + A^{-2})D).$$
\end{enumerate}
\end{definition}

\begin{lemma}\label{7.3}\ \\
Let $B^d(F)$ denote the set of link diagrams in $F$ without vertices 
and without trivial components (but allowing $\emptyset$). We can
identify $B^d(F)$ with the set $B(F)$ with framing ignored.
$B^d(F)$ is a subset of the set of link diagrams, so we have a homomorphism 
$\phi: RB^d(F) \to {\cal S}{\cal D}(F;R,A)$ defined by
associating to a link, $\gamma \in B^d(F)$, in $F$ its class in 
${\cal S}{\cal D}(F;R,A)$.\\
Then $\phi$ is an isomorphism.
\end{lemma}

\begin{proof} For any $D \in \cal D$ we can use the first relation to 
eliminate all crossings, and the second to eliminate trivial components of
$D$. Thus $\phi$ is an epimorphism.\\
To show that it is a monomorphism we will construct the inverse 
map, $\psi$.\\
First we define a map $\hat\psi: R{\cal D} \to RB^d(F)$.
Let $D\in \cal D$. We define $\hat\psi(D)$ as follows:\\
Choose any ordering $p_1,...,p_n$ of crossings of $D$, and use the
formula $D^{p_i}_{{\psfig{figure=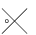}}}=
AD^{p_i}_{{\psfig{figure=L0nmaly.eps}}} + 
A^{-1} D^{p_i}_{{\psfig{figure=Linftynmaly.eps}}}$, for each crossing,
until all crossings are eliminated.
The upper index denotes the crossing
at which we perform a smoothing (crossing elimination). 
The result does not depend on the order of the crossings since we can make
any transposition of adjacent (with respect to ordering) pairs 
and get the same result:\\
\ \\
$(D^p_{{\psfig{figure=Lmarkednmaly.eps}}})^{q}_{{\psfig{figure=Lmarkednmaly.eps}}} =
A(D^p_{{\psfig{figure=L0nmaly.eps}}})^q_{{\psfig{figure=Lmarkednmaly.eps}}} +
A^{-1}(D^p_{{\psfig{figure=Linftynmaly.eps}}})^q_{{\psfig{figure=Lmarkednmaly.eps}}}
=\\
\ \\
A^2(D^p_{{\psfig{figure=L0nmaly.eps}}})^q_{{\psfig{figure=L0nmaly.eps}}} +
(D^p_{{\psfig{figure=L0nmaly.eps}}})^q_{{\psfig{figure=Linftynmaly.eps}}} +
(D^p_{{\psfig{figure=Linftynmaly.eps}}})^q_{{\psfig{figure=L0nmaly.eps}}} +
A^{-2}(D^p_{{\psfig{figure=Linftynmaly.eps}}})^q_
{{\psfig{figure=Linftynmaly.eps}}}$\\ 
\ \\
and \\
$(D^q_{{\psfig{figure=Lmarkednmaly.eps}}})^{p}_{{\psfig{figure=Lmarkednmaly.eps}}} =   
A(D^q_{{\psfig{figure=L0nmaly.eps}}})^p_{{\psfig{figure=Lmarkednmaly.eps}}} + 
A^{-1}(D^q_{{\psfig{figure=Linftynmaly.eps}}})^p_{{\psfig{figure=Lmarkednmaly.eps}}}  
=\\
\ \\                                                                      
A^2(D^q_{{\psfig{figure=L0nmaly.eps}}})^p_{{\psfig{figure=L0nmaly.eps}}} + 
(D^q_{{\psfig{figure=L0nmaly.eps}}})^p_{{\psfig{figure=Linftynmaly.eps}}} + 
(D^q_{{\psfig{figure=Linftynmaly.eps}}})^p_{{\psfig{figure=L0nmaly.eps}}} + 
A^{-2}(D^q_{{\psfig{figure=Linftynmaly.eps}}})^p_ 
{{\psfig{figure=Linftynmaly.eps}}}$\\ 
\ \\
After smoothing all crossings we eliminate trivial components by the relation
 $D\sqcup T_1 = (-A^2 -A^{-2})D$ (there is no ambiguity in the reduction).
Thus $D$ is uniquely expressed as a linear combination of elements of
$B^d(F)$, and we define $\hat\psi(D)$ as this linear combination (which 
lies in $RB^d(F)$). Therefore $\hat\psi$ is well defined.
Now $\hat\psi$ descends to $\psi: {\cal S}{\cal D}(F;R,A) \to RB^d(F)$ because
$\hat\psi ( \parbox{0.5cm}{\psfig{figure=Lmarkednmaly.eps}} - 
A\parbox{0.5cm}{\psfig{figure=L0nmaly.eps}} - 
A^{-1} \parbox{0.5cm}{\psfig{figure=Linftynmaly.eps}})=0$
and $\hat\psi ( D\sqcup T_1 +(A^2 + A^{-2})D)=0$. 
Now, obviously, $\psi\phi = Id$, thus $\phi$ is a monomorphism.
\end{proof}
  Our goal is to prove that $B(F)$ is a basis of the Kauffman bracket
skein module  ${\cal S}_{2,\infty}(M;R,A)$, where $M$ is an oriented
3-manifold which is an $I$-bundle over a surface $F$. Because we would like to
consider the case of orientable and unorientable surface simultaneously,
it is convenient to consider half-integer framings of links, that
is, to allow embedded M\"obius bands. This suggests the following
definition.
\begin{definition}\label{IX.7.4}
Let $M$ be any oriented 3-manifold, $\bar{\cal L}_{fr}$ the set of embeddings
of annuli and M\"obius bands in $M$ (up to an ambient isotopy of $M$) 
and $\bar R$ 
a commutative ring with identity with a chosen invertible element $\bar A$ 
(we define $A=-\bar A^2$ and we will often write $\sqrt{-A}$ for $\bar A$). 
Let $\bar R \bar{\cal L}^{fr}$ denote a free $\bar R$ module over 
$\bar{\cal L}^{fr}$ and let $\bar S_{2,\infty}$ denote the submodule of
$\bar R\bar{\cal L}^{fr}$ generated by expressions 
$L_+ -AL_- - A^{-1}L_{\infty}$,
and $L^{1/2} - (\sqrt{-A})^3L$, where $L^{1/2}$ denotes $L$ with its 
framing twisted by a half twist in a positive direction. 
As before, for convenience, we allow the empty knot, $\emptyset$, and
add the relation $T_1=(-A^2-A^{-2})\emptyset$. \\
Then we define $\bar{\cal S}_{2,\infty}(M,\bar R,\bar A)= 
\bar R \bar{\cal L}_{fr}/\bar S_{2,\infty}.$
\end{definition}

Consider the $\bar R$-homomorphism
$g: {\cal S}{\cal D}(F;\bar R, A) \to \bar{\cal S}_{2,\infty}(M,\bar R,\bar A)$
defined on the basic elements $\gamma \in B^d(F)$ by $g(\gamma)=
\gamma^{fr}$ where $\gamma^{fr}$ is a framed link obtained from
$\gamma$ by giving it the blackboard framing (it may be an annulus or a 
M\"obius band). Using our skein relations,
in a similar manner as before, we see that $g$ is an epimorphism.
If $D$ is any marked diagram we can describe the framed link $g(D)$ as
follows: we resolve every crossing of $D$ according to the rule given
in Fig. 7.1 and giving the link $g(D)$ the blackboard framing
(the orientation of $M$ in a neighborhood of the crossing should
agree with that of $R^3$ from Fig. 7.1).\ \\
\ \\
\centerline{\psfig{figure=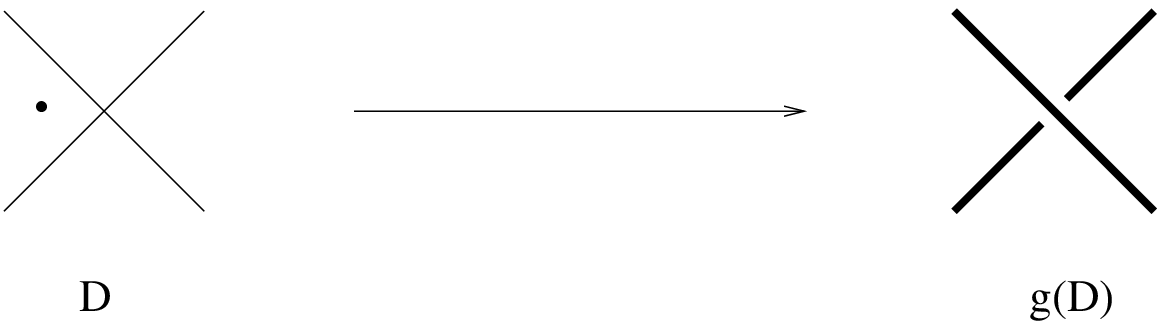}}
\begin{center}
Fig. 7.1.
\end{center}

Consider now the following lemma
concerning Reidemeister moves on diagrams.
\begin{lemma}\label{7.5}
Consider the following moves $\bar R_1, \bar R_2, \bar R_3$ on marked diagrams.
In ${\cal S}{\cal D}(F;\bar R, A)$ they satisfy:
\begin{enumerate}
\item[($\bar R_1$)]\ {\psfig{figure=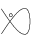}} 
$= - A^3$ {\psfig{figure=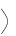}}
and {\psfig{figure=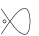}} 
$=- A^{-3}$ {\psfig{figure=Lline.eps}}, where
$\bar R_1$({\psfig{figure=Lline.eps}})$=$ {\psfig{figure=Ltwistgora.eps}} or
{\psfig{figure=Ltwistdol.eps}}. 
\item[($\bar R_2$)]\ $\bar R_2(D)=D$, where 
$\bar R_2$({\psfig{figure=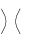}})
$=$
{\psfig{figure=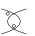}}.
\item[($\bar R_3$)]\  $\bar R_3(D)=D$, where 
$\bar R_3$({\psfig{figure=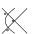}})
$=$
{\psfig{figure=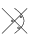}}.
\end{enumerate} 
\end{lemma}
\begin{proof}
\begin{enumerate} 
\item[($\bar R_1$)]\ \\ 
{\psfig{figure=Ltwistgora.eps}} $= A$\ {\psfig{figure=Lline.eps}}\ $\sqcup O 
+A^{-1}${\psfig{figure=Lline.eps}} $(A(-A^2-A^{-2}) + A^{-1}$
{\psfig{figure=Lline.eps}} $= - A^3$ {\psfig{figure=Lline.eps}}.\\
{\psfig{figure=Ltwistdol.eps}} $= A^{-1}$ {\psfig{figure=Lline.eps}}$\sqcup O
A${\psfig{figure=Lline.eps}} $(A+A^{-1}(-A^2-A^{-2})$ 
$=- A^{-3}$ {\psfig{figure=Lline.eps}}.
\item[($\bar R_2$)]\ \\
{\psfig{figure=R2marked.eps,height=0.4cm}}
$=A$ {\psfig{figure=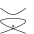,height=0.4cm}} $+ A^{-1}$
{\psfig{figure=L-nmaly.eps,height=0.4cm}} 
$= A(A$ {\psfig{figure=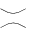,height=0.4cm}} $+ A^{-1}$
{\psfig{figure=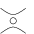,height=0.4cm}}) 
$+A^{-1}(A$ {\psfig{figure=L2lines.eps,height=0.4cm}}
$+A^{-1}$ {\psfig{figure=L2linesh.eps,height=0.4cm}}) 
$= (A^2 + AA^{-1}(-A^2 -A^{-2})+A^{-2})$ 
{\psfig{figure=L2linesh.eps,height=0.4cm}} $+$
{\psfig{figure=L2lines.eps,height=0.4cm}}
$=${\psfig{figure=L2lines.eps,height=0.4cm}}.
\item[($\bar R_3$)]\ \\
{\psfig{figure=R3markedb.eps,height=0.4cm}}$= A$
{\psfig{figure=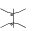,height=0.4cm}}
$+ A^{-1}$
{\psfig{figure=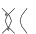,height=0.4cm}} 
$=A$
{\psfig{figure=R3markedbAsplit.eps,height=0.4cm}}
$+ A^{-1}$
{\psfig{figure=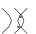,height=0.4cm}}
$=$ 
{\psfig{figure=R3markede.eps,height=0.4cm}}. 
We use here the invariance under $\bar R_2$
moves.
\end{enumerate}
\end{proof}
To use the lemma in the proof that $g$ is a monomorphism, we need a
variant of Reidemeister's theorem for marked diagrams:

\begin{proposition}\label{7.6}
Let $\hat g: {\cal D} \to {\cal L}_{fr} $ be a map given by Fig. 5.1.
Then two marked diagrams, $D_1$ and $D_2$, represent 
the same framed link, $\hat g(D_1)=\hat g(D_2)$, \\ if and only if\\
one can go from $D_1$ to $D_2$ using Reidemeister moves $\bar R_i^{\pm 1}$ 
and an isotopy of $F$, and additionally, for corresponding link components
of $D_1$ and $D_2$, their Tait numbers are the same. One should notice here
that for a knot diagram the Tait number is independent on orientation
of the knot. Precisely for a knot diagram $D$ we define 
Tait($D$)$= \Sigma_p sgn(p)$, where 
sgn({\psfig{figure=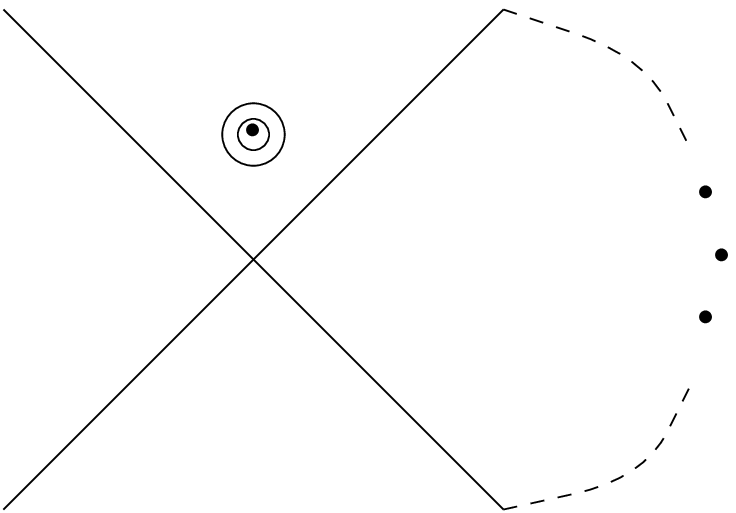,height=0.4cm}})$=1$ and 
sgn({\psfig{figure=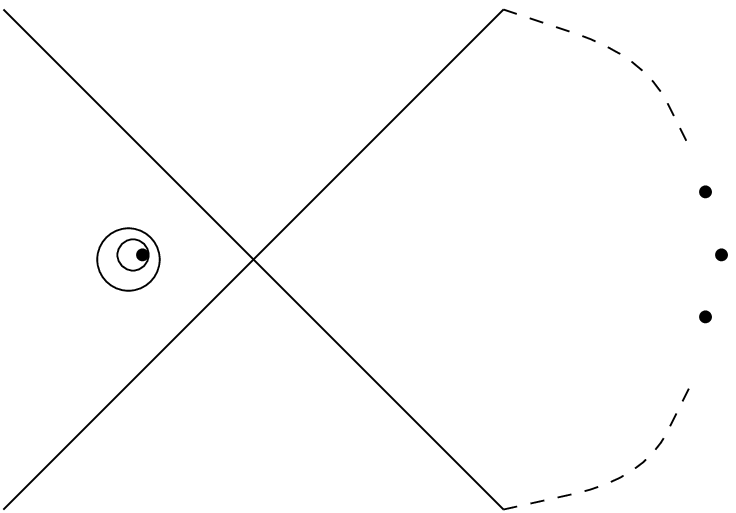,height=0.4cm}})$=-1$.

\end{proposition}

\begin{proof} The proposition can be deduced from the classical Reidemeister
theorem and the result from the PL topology; Theorem 6.2 in
\cite{Hud}\footnote{
It follows from the theorem that if $C$ is a compact subset of a manifold 
$M$ and $F:M\times I \to M$ is the isotopy of $M$ then there is another
isotopy $\hat F:M\times I \to M$ such that\\
 $F_0 ={\hat F}_0$,
$F_1/C = {\hat F}_1/C$ and there exists a number $N$ such that
the set\\ $\{x\in M \ |\ \hat F/\{x\}\times (k/N,(k+1)/N) \ is \ 
not \ constant\}$
sits in a ball embedded in $M$.}
\end{proof}
Our goal is to show that the epimorphism
$g: {\cal S}{\cal D}(F;\bar R, A) \to \bar{\cal S}_{2,\infty}(M,\bar R,\bar A)$
is an isomorphism.
We use Lemma 7.5 and Proposition 7.6 to construct the map inverse to
$g$.\
Let $\hat h: \bar R \bar{\cal L}_{fr} \to {\cal S}{\cal D}(F;\bar R,A)$ be a 
homomorphism defined as follows: choose a representative of a link 
$L\in \bar{\cal L}_{fr}$ which has a regular projection on $F$. Let $D_L$ be
a marked diagram on $F$ constructed as in Fig. 5.1, and let $t(L)$ be
the number (possibly half-integer) of positive twists which should 
be performed on the blackboard
framing of $D_L$ to get the framing of $L$. Then we define
$\hat h(L)=(-A^3)^{t(L)}D_L$. $\hat h(L)$ is well defined by Lemma 7.5 
and Proposition 7.6. 
Furthermore, $\hat h(L_+ -AL_- -A^{-1}L_{\infty})=0$, 
$\hat h(L \sqcup T_1 +(A^2 + A^{-2})L)=0$ and 
$\hat h(L^{1/2}-{\sqrt {-A^3}}L)=0$ so $\hat h$ descends to
$h:{\cal S}_{2,\infty}(M;R,A) \to {\cal S}{\cal D}(F;\bar R,A)$.
Of course $hg=Id$ so $g$ is a monomorphism, as required.
\end{proof}
We can now complete the proof of Theorem 5.1. Because 
$g: {\cal S}{\cal D}(F;\bar R, A) \to \bar{\cal S}_{2,\infty}(M,\bar R,\bar A)$
is an isomorphism, therefore $g(B^d(F))$ is a basis of 
$\bar{\cal S}_{2,\infty}(M,\bar R,\bar A)$. On the other hand $B(F)$, whose
elements may differ from elements of $g(B^d(F))$ only by framing, also forms
a basis of $\bar{\cal S}_{2,\infty}(M;\bar R,\bar A)$. Thus they are linearly
independent in ${\cal S}_{2,\infty}(M;R,A)$. Because $B(F)$ generates
${\cal S}_{2,\infty}(M;R,A)$ it is a basis of this module. The proof of
Theorem 7.1 is completed.

As an immediate corollary of Theorem 7.1 we obtain the structure of
 KBSM of the projective space $RP^3$. This result was also obtained 
independently by J.~Drobotukhina \cite{Dr}. 
\begin{corollary}\label{IX.7.7} 
${\cal S}_{2,\infty}(RP^3;R,A)=R\oplus R$. As a basis
of KBSM we can take $\emptyset$ and a generator of the fundamental group
of $RP^3$.
\end{corollary}
\begin{proof}
By Proposition 6.2(i) ${\cal S}_{2,\infty}(RP^3;R,A)=
{\cal S}_{2,\infty}(RP^3-int(D^3);R,A)$ and $RP^3-int(D^3)$ is equal to 
the twisted $I$-bundle over a projective plane ($RP^2\hat I$. By
Theorem 5.1, ${\cal S}_{2,\infty}(RP^2\hat I;R,A)$ is a free $R$-module
with basis $B(RP^2)$, which has two elements: the empty knot and the
noncontractible curve on $RP^2$. 
\end{proof}
One can generalize Theorem 7.1 to relative skein modules, as long
as boundary points of relative links are on the same boundary
component of a manifold.

\begin{definition}[Relative Kauffman Bracket Skein Module]\label{IX.7.8}\ \\
Let $x_1,x_2,...,x_{2n}$ be a set of $2n$ (framed) points in $\partial M$,
where $M$ is an oriented 3-manifold. Let ${\cal L}_{fr}(n)$ be a family
of relative framed links in $(M,\partial M)$ such that $L\hat \partial M =
\partial L =\{x_i\}$, considered up to an ambient isotopy fixing $\partial M$.
Let $R$ be a commutative ring with identity and $A$ its 
invertible element. Let $S_{2,\infty}(n)$ be 
the submodule of $R{\cal L}_{fr}(n)$ generated by the Kauffman bracket
skein relations. We define the Relative Kauffman Bracket Skein Module (RKBSM)
as the quotient:
$${\cal S}_{2,\infty}(M,\{x_i\}_1^{2n};R,A)= 
R{\cal L}_{fr}(n)/S_{2,\infty}(n)$$

\end{definition}

We list below a few useful properties of relative skein modules:
\begin{proposition}\label{7.9}
\begin{enumerate}
\item[(a)] There is a functor from the category of
oriented 3-manifolds with $2n$ framed points on the boundary and 
orientation preserving embeddings (up to ambient isotopy fixed on the boundary)
to the category of $R$-modules (with a specified invertible element $A \in R$).
The functor sends an embedding $i: (M,\{x_i\}_{i=1}^{2n}) \to
(N,\{y_i\}_{i=1}^{2n})$ into $R$-modules morphism
${\cal S}_{2,\infty}(M,\{x_i\}_1^{2n};R,A) \to
{\cal S}_{2,\infty}(N,\{y_i\}_1^{2n};R,A)$.
\item[(b)] Adding a 3-handle to $M$ (outside $x_i$) does not change the RKBSM,
and adding a 2-handle is adding only relations to RKBSM (handle slidings 
yield relations); compare Proposition 6.2(2). 
\item[(c)] The relative KBSM depends only on the distribution of boundary
points $\{x_i\}$ among boundary components of $M$, but not on the exact
position of $\{x_i\}$. In particular if $\partial M$ is connected, we can
write shortly ${\cal S}_{2,\infty}(M,n;R,A)$ instead of 
${\cal S}_{2,\infty}(M,\{x_i\}_1^{2n};R,A)$
\item[(d)] The relative KBSM satisfies the Universal Coefficient Property,
compare Proposition 6.2(4).
\item[(e)] For a disjoint sum of 3-manifolds we have:
$${\cal S}_{2,\infty}(M_1\sqcup M_2,\{x_i,y_i\}_1^{2n};R,A)= 
{\cal S}_{2,\infty}(M_1,\{x_i\}_1^{2n};R,A) 
\otimes {\cal S}_{2,\infty}(M_2,\{y_i\}_1^{2n};R,A).$$
\end{enumerate}
\end{proposition}

\begin{theorem}\label{7.10} Let $M= F\bar{\times} I$ that is
 $M= F\times I$ or $M= F\hat{\times} I$, then
\begin{enumerate}
\item [(a)] Let $\partial F \neq \emptyset$ then
${\cal S}_{2,\infty}(M,\{x_i\}_1^{2n};R,A)$ is a free $R$-module. 
Consider all $x_i$ to
lie on $\partial F \times \{\frac{1}{2}\}$ then the basis of the
module ${\cal S}_{2,\infty}(M,\{x_i\}_1^{2n};R,A)$ is composed of
relative links on $F$ without trivial components.
\item [(b)] In the case of $F_{g,0}$ closed surface of genus $g$
 ($F\neq S^2$) the situation is more
delicate so we stop on the following observation:\\ 
${\cal S}_{2,\infty}(F_{g,0}\hat{\times} I;\{x_i\}_1^{2n};R,A)=
{\cal S}_{2,\infty}(F_{g,1}\hat{\times} I;\{x_i\}_1^{2n};R,A)/(I)$ where\\ 
$F_{g,1}=F_{g,0}-int(D^2)$ and assuming $x_i \in \partial D^2$, ideal $(I)$
is generated by moves in which arcs go above $D^2$.
\end {enumerate}
\end{theorem}
\begin{proof} The proof of (a) is the same as that of Theorem 7.1;
as before relative link diagrams representing the same link are
related by Reidemeister moves. In the case (b) it is no longer true as
we need also handle sliding. $F_{g,0}\hat{\times} I$ is obtained
from $F_{g,1}\hat{\times} I$ by adding the 2-handle along $\partial D^2$.
Now (b) follows from Proposition 7.9(b).  
\end{proof}
In the case $F$ is a closed surface the question whether 
${\cal S}_{2,\infty}(F \hat \times I ,\{x_i\}_1^{2n};R,A)$ is free is 
open in general. If not all $x_i$ lie on the same boundary component
of $F \times I $ the the skein module has a torsion in the case
of $F$ being a sphere or a torus; compare Section 6.
 We propose the following conjecture,
which we are able to confirm only for $F$ being a torus and $n=1$.
\begin{conjecture}\label{IX.7.11}
Let $F$ be a closed surface and $x_i\in F\times \{0\}$ for any $i$, then
the skein module ${\cal S}_{2,\infty}(F \times I ,\{x_i\}_1^{2n};R,A)$ 
is free.
\end{conjecture}

\begin{corollary}\label{7.12}\ 
\begin{enumerate} 
\item [(a)] ${\cal S}_{2,\infty}(D^2\times I,\{x_i\}_1^{2n};R,A)$ is
a free $R$ module of $\frac{1}{n+1} {{2n}\choose{n}}$ free generators.

\item [(b)] 
${\cal S}_{2,\infty}(annulus\ \times I,\{x_i\}_1^{2n};R,A)$ is
a free $R[\alpha]$ module with ${2n}\choose{n}$ free generators,
where $\alpha$ is represented by a longitude of the annulus.

\end{enumerate}
\end{corollary}
\begin{proof}
Corollary 7.12 (a) describes, well known, module structure of the 
Temperley-Lieb algebra (basis having Catalan number of generators).
(b) follows from work of Jones and Tom Dieck. We give here a self-sufficient
proof. In lieu of Theorem 7.10,
it suffices to count the crossless connections of $2n$ points in the disc
and annulus. We offer here an amazingly simple calculation for both cases
simultaneously. Let $C_n$ be the number of connections in the disc and
$D_n$ the number of connections in the annulus (all points $x_i$ are
on the ``outside" circle of the annulus.
Connections in the disc cut it into $n+1$ pieces; to get connection in
the annulus we have to put a ``table" (remove a disk) from $D^2$; thus
$D_n = (n+1)C_n$. On the other hand any arc of a connection in the annulus 
has the first point (with respect to fixed orientation of the annulus),
and any choice of $n$ points leads to unique connection, for which given
points are first. Therefore $D_n= {{2n}\choose {n}}$ and thus $C_n=
\frac{1}{n+1} {{2n}\choose{n}}$.
\end{proof}
We finish this section by offering the following very useful observation
(compare \cite{P-S-2}).

\begin{proposition}\label{7.13}
Consider a 3-manifold $(M,\{x_i\}_{i=1}^{2n})$ and 
let $x_{2n+1}$ and $x_{2n+2}$ lie on the same boundary component of $M$.
Consider the $R$-homomorphism of RKBSM
$$i_{\#}: {\cal S}_{2,\infty}(M,\{x_i\}_{i=1}^{2n},R,A) \to
{\cal S}_{2,\infty}(M,\{x_i\}_{i=1}^{2n+2},R,A)$$
 generated by the identity map and with convention that $i_{\#}(L)$
have $x_{2n+1}$ connected to $x_{2n+2}$ by a framed arc close to
boundary (we push out of the boundary framed arc joining $x_{2n+1}$ and
$x_{2n+2}$ in $\partial M$).\ \
Then\\
$i_{\#}$ is a monomorphism if one assumes that
$A^2+A^{-2}$ is not an annihilator of any non-zero element of 
${\cal S}_{2,\infty}(M,\{x_i\}_{i=1}^{2n},R,A)$ 
(i.e. $(A^2+A^{-2})x=0 \Rightarrow x=0$
\end{proposition}
\begin{proof} Consider the $R$-homomorphism 
$$i_{\#}': {\cal S}_{2,\infty}(M,\{x_i\}_{i=1}^{2n+2},R,A) \to
{\cal S}_{2,\infty}(M,\{x_i\}_{i=1}^{2n},R,A)$$ given by connecting
$x_{2n+1}$ and $x_{2n+2}$ in $\partial M$ and pushing it inside $M$.
Now clearly $i_{\#}'i_{\#}(L)=(-A^2-A^{-2})(L)$, thus 
$i_{\#}'i_{\#}(u)=(-A^2-A^{-2})(u)$ for any $u\in 
{\cal S}_{2,\infty}(M,\{x_i\}_{i=1}^{2n},R,A)$. Therefore $i_{\#}'i_{\#}$
is a monomorphism iff $A^2+A^{-2}$ is not an annihilator of any non-zero 
element of ${\cal S}_{2,\infty}(M,\{x_i\}_{i=1}^{2n},R,A)$. If 
$i_{\#}'i_{\#}$ is a monomorphism then $i_{\#}$ is a monomorphism.
\end{proof}
  
\section{Torsion in KBSM}\label{IX.8}
In all of the  examples above the module is torsion free except in 
the case of $S^1 \times S^2$.
In fact a non-separating $S^2$ in $M$ always yields a torsion in  
${\cal S}_{2,\infty}(M)$. It is enough to use the framing relation to see
a torsion: Let $L$ be a framed link cutting a non-separating $S^2$ exactly
in one point. We can twist $S^2$ twice, twisting also the framing of $L$
twice and then undo this by an isotopy of $M$. Thus $(A^6-1)L=0$ in
${\cal S}_{2,\infty}(M)$. It is less obvious that a separating $S^2$ can
often yield a torsion.
\begin{conjecture}[\cite{Kir}] \label{8.1}
If $M=M_1\# M_2$, where $M_i$ is not equal to $S^3$, possibly with holes,
then ${\cal S}_{2,\infty}(M)$ has a torsion element.
\end{conjecture}

We are able to prove Conjecture 8.1 only partially.
\begin{theorem}\label{8.2}
If $M_1$ and $M_2$ have first homology groups that are not 
2-torsion groups, then the conjecture holds.
\end{theorem}

The proof will be divided into two lemmas, each of which is of 
considerable interest.
\begin{lemma}\label{6.3}
Consider the $Z$-epimorphism $\phi$ from $R{\cal L}^{fr}$ onto $ZH_1(M,Z)$
given by $\phi (L) =(-1)^{com(L)}\sum_{L_i\in or(L)}|L_i|$ where the sum
is taken over all possible orientations of $L$ and $|L_i|$denotes the 
homology class of the oriented link $L_i$. Furthermore $\phi (Aw)=-\phi (w)$
(i.e. $A \to -1$). Then $\phi$ descends to the map $\hat\phi:
{\cal S}_{2,\infty}(M) \to ZH_1(M,Z)$. 
More succinctly, although imprecisely, one can say that the Kauffman 
bracket skein module is more delicate than the homology up to
orientation.
\end{lemma}

\begin{proof}
It suffices to show that $\phi$ sends skein expressions to zero.
Because $A$ is sent to $-1$ then $\phi (L^{(1)}+A^3L)=0$.
Furthermore $\phi (L_+ -AL_0-A^{-1}L_{\infty})=\phi (L_+ +L_0 +L_{\infty})$
and among three links $L_+, L_0, L_{\infty}$ exactly one has more
components than the other two (hence twice as many orientations). 
Comparing these orientations, we see that all expressions reduce
to zero (in $ZH_1(M,Z)$).
\end{proof}

\begin{lemma}\label{8.4}
\begin{enumerate}
\item [(a)] Let $(M,K)=(M_1,K_1)\# (M_2,K_2)$ and \\
$(M,K')=(M_1,K_1) (\#,\sqcup) (M_2,K_2)$, where ($\#,\sqcup$) means that
we consider the disjoint sum of knots $K_1 \sqcup K_2$ in the connected sum
$M_1 \# M_2$. Then  $(A^4-1)(\mu K -K')=0$ in the skein module 
${\cal S}_{2,\infty}(M)$, where $\mu = -A^2-A^{-2}$.
\item [(b)] If $2|K_i|\neq 0$ in $H_1(M_i,Z)$, $i=1,2$, then $\mu K -K'
\neq 0$ in ${\cal S}_{2,\infty}(M)$.
\end{enumerate}
\end{lemma}
\ \\
\ \\
\centerline{\psfig{figure=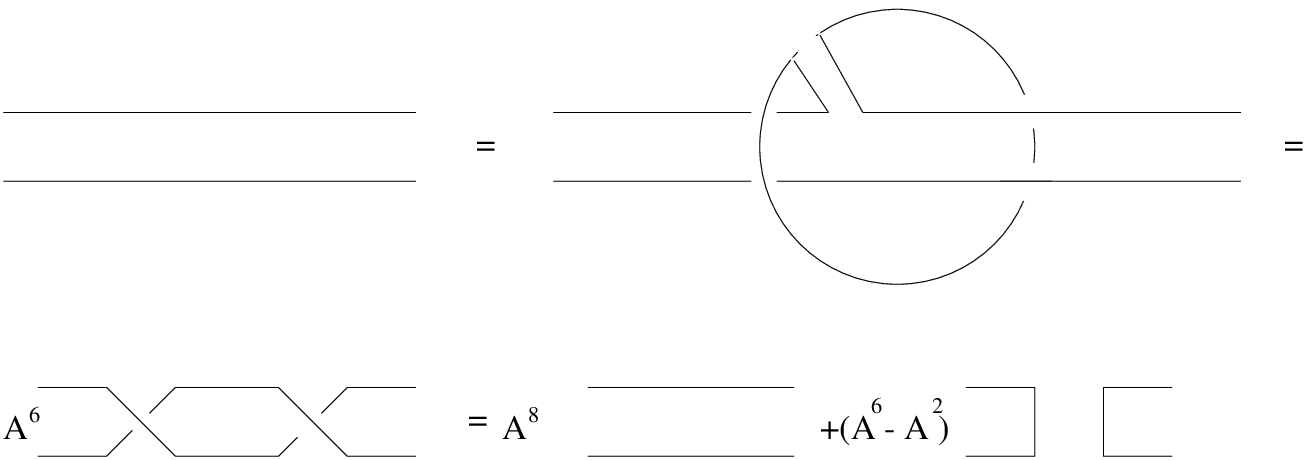}}
\begin{center}
Fig. 8.1.
\end{center}

\begin{proof}
\begin{enumerate}
\item [(a)] 
By using the ``second" side of $S^2$ one gets the identity from 
Fig. 5.2. After reducing the right-hand-side of the equation, using skein
relations, one gets in ${\cal S}_{2,\infty}(M)$ the identity:
$ K= A^8K + (A^6-A^2)K'$ and therefore $0=(A^8-1)K+(A^6-A^2)K'=
A^2(A^4-1)((A^2+A^{-2})K+K')=-A^2(A^4-1)(\mu K -K')$. Therefore for 
$f_2=\mu K -K'$ we have $(A^4-1)f_2=0$ in ${\cal S}_{2,\infty}(M)$, 
as required\footnote{ $f_2$ is related to Jones-Wenzl second idempotent in the 
Temperley-Lieb algebra.}.

\item [(b)]
Consider the epimorphism $\hat\phi$ from Lemma x.8.
$\hat\phi(\mu K -K')= -2 \phi (K) - \phi (K')$ in $ZH_1(M,Z)$.
If $K_i^{\pm}$ denotes two possible orientations of 
$K_i$ and $K^+=K^+_1 \# K^+_2$,
and $K^-=K^-_1 \# K^-_2$, then 
$\hat\phi(\mu K -K')= -|K_1^+ \sqcup K_2^+|+|K_1^+ \sqcup K_2^-|
-|K_1^- \sqcup K_2^-|+|K_1^- \sqcup K_2^+|$, which is equal to zero
if and only if $|K_1^+|=|K_1^-|$ or $|K_2^+|=|K_2^-|$.
\end{enumerate}
\end{proof}

\begin{example}\label{6.5}
\begin{enumerate}
\item[(a)] ${\cal S}_{2,\infty}(S^1\times S^1\times S^1)$ has a torsion.
\item[(b)] The Kauffman bracket skein module of the 
double of the complement of the figure eight knot (Listing knot) has a torsion.
\end{enumerate}
\end{example}
To prove (a) the idea of ``symmetric homologies" described in Lemma 6.3
is crucial (``symmetric homology" is the image of $\hat\phi$ in 
$ZH_1(M,Z)$). For (b) the new idea is the observation by Bullock that
there is an algebra homomorphism from the skein algebra ($A=-1$) to $C$
yielded by a homomorphism of the fundamental group of the manifold to
$SL(2,C)$. Further we use the existence of hyperbolic structure on the
complement of the Listing knot (it is probably the first connection
between Jones type invariants and hyperbolic structures).
  
\section{Kauffman bracket skein algebra of $F\times I$ and skein algebra
of a group}
\markboth{\hfil{\sc Skein modules}\hfil}
{\hfil{\sc Skein algebra of a group}\hfil}

We discuss in this section a possible algebra structure for a Kauffman
bracket skein module. We can get the structure either by considering 
a special $M$ (e.g. $F\times I$ or $RP^3$) or a special $A$ in $R$
(e.g. $A=\pm 1$). We allow the empty knot, $\emptyset$, in order to 
have a unit of an algebra.

\begin{theorem}{[Bullock]}\label{9.1}
If ${\cal S}_{2,\infty}(F\times [0,1])$ is given the multiplication 
$L_1\cdot L_2$ defined by placing $L_1$  above $L_2$ 
then the resulting algebra is finitely generated \cite{Bu-3} and the
minimal number of generators $r({\cal S}_{2,\infty}(F\times [0,1])$ is
no more than $2^{rank(H_1(F))}-1$. 
\end{theorem}
More precisely we have:
\begin{theorem}[(\cite{P-S-2})]\label{9.2}\ \\ 
\begin{enumerate} 
\item [(i)] 
$r({\cal S}_{2,\infty}(F\times I;Z[A^{\pm 1}],A))=2^{rank(H_1(F))}-1$ 
\item[(ii)] If $R=Z[A^{\pm 1}, (A^2 + A^{-2})^{-1}]$ then 
$$r({\cal S}_{2,\infty}(F\times I;R,A))= d+ {d\choose 2} + {d\choose 3},$$ 
where $F$ is a disk with $d$ holes. 
\end{enumerate} 
\end{theorem} 
\begin{theorem}
${\cal S}_{2,\infty}(L(2,1);R,A)$ has an 
algebra structure and as an algebra it is isomorphic to 
$R[\alpha]/(\alpha^2-A^3\frac{ 
A^4-A^{-4}}{A-A^{-1}})$
\end{theorem}
\begin{theorem}
\begin{enumerate} 
\item [(i)] 
${\cal S}_{2,\infty}(F\times I;R,A)$ has no zero divisors, provided $R$ 
has no zero divisors.
\item [(ii)] The center of the algebra ${\cal S}_{2,\infty}(F\times I;R,A)$ 
is a subalgebra generated by the boundary components of $F$. 
\item [(iii)] 
If $M$ is a twisted $I$ bundle over an unorientable surface $F$ 
(having an even number of projective planes as factors) then 
${\cal S}_{2,\infty}(M;R,-1)$ has no zero divisors, provided $R$ 
has no zero divisors and $F$ is not a Klein bottle.
\item [(iv)] ${\cal S}_{2,\infty}(T^2\times I;Z,-1)$ 
is a unique factorization domain. 
\item [(v)]
The skein algebra ${\cal S}_{2,\infty}(F\times I;C,-1)$, is isomorphic 
to the coordinate ring of the $SL(2,C)$ character variety of 
the fundamental group of the surface (it solves Bullock conjecture for
$F\times I$). .
\end{enumerate} 
\end{theorem}

The exact structure of the Kauffman bracket skein algebra is computed
only for ``small surfaces" $F_{g,d}$ where $g$ is the genus of the surface and
$d$ the number of its boundary components \cite{B-P}. 

\begin{proposition} 
\begin{enumerate} 
\item [(1)]
${\cal S}_{2,\infty}(F_{0,1}\times I;R,A)=R$
\item [(2)]
${\cal S}_{2,\infty}(F_{0,2}\times I;R,A)=R[z]$
\item [(3)] ${\cal S}_{2,\infty}(F_{0,3}\times I;R,A)=R[x,y,z]$
\end{enumerate}
\end{proposition}
For surfaces which contain simple closed curves not parallel to the boundary
the skein algebra is not commutative (see Theorem 5.4(ii)).\\ 
In what follows, use notation $<x,y,...>$ to denote noncommutative variables.

\begin{theorem}\ \\
\begin{enumerate}  
\item [(1)]
$${\cal S}_{2,\infty}(F_{1,1}\times I);R,A)= R<x,y,z>/I_{1,1}$$ where 
$I_{1,1}$ is the ideal generated by ``$A$-commutators"\\
$Axy - A^{-1}yx\ -\ (A^2-A^{-2})z$ \\ 
$Ayz - A^{-1}zy\ -\ (A^2-A^{-2})x$ \\ 
$Azx - A^{-1}xz\ -\ (A^2-A^{-2})y$ \\ 
\item [(2)]
${\cal S}_{2,\infty}(F_{1,0}\times I);R,A)=R<x,y,z>/I_{1,0}$, \ where 
$I_{1,0}$ is the ideal generated by $I_{1,1}$ and the long relation
$$A^2x^2 + A^{-2}y^2 + A^2z^2-Axyz -2(A^2+A^{-2}).$$ 
\end{enumerate}
\end{theorem}
\begin{theorem}\label{x.x}\ \\ 
${\cal S}_{2,\infty}(F_{0,4}\times I;R,A)= 
R[a_1,a_2,a_3,a_4]<x_1,x_2,x_3>/(I_{0,4})$, 
where  $I_{0,4}$ is an ideal generated by ``$A^2$ bracket" relations: 
$$A^2x_ix_{i+1} -A^{-2}x_{i+1}x_i =(A^4-A^{-4})x_{i+2} + 
(A^2-A^{-2})(a_ia_{i+1} + a_{i+2}a_4),$$ 
where $i=1,2, 3$ and indices are taken modulo $3$, 
and by the ``long relation": 
$$A^4x_1^2 + A^{-4}x_2^2 + A^4x_3^2 -A^2x_1x_2x_3 + A^2x_1(a_2a_3+a_1a_4) + 
A^{-2}x_2(a_1a_3+a_2a_4) +$$ 
$$A^{-2}x_3(a_1a_2+a_3a_4)+ a_1^2 +a_2^2 +a_3^2 + 
a_4^2 + a_1a_2a_3a_4 -A^4-A^{-4} -2$$
\end{theorem}
\begin{theorem} 
$${\cal S}_{2,\infty}(F_{1,2}\times I;R,A)= R[a]<x,y,z,x',y',z'>/ 
I_{1,2}$$ where $I_{1,2}$ is the ideal generated by $A$-commutation 
relations:\\ 
$Axy-A^{-1}yx = (A^2-A^{-2})z$, 
$Ayz-A^{-1}zy = (A^2-A^{-2})x$,\\ 
$Azx-A^{-1}xz = (A^2-A^{-2})y$, 
$Axy'-A^{-1}y'x = (A^2-A^{-2})z'$,\\ 
$Ax'y-A^{-1}yx' = (A^2-A^{-2})z'$ 
$Ayz'-A^{-1}z'y = (A^2-A^{-2})x'$,\\ 
$Ay'z-A^{-1}zy' = (A^2-A^{-2})x'$, 
$Azx'-A^{-1}x'z = (A^2-A^{-2})y'$,\\ 
$Az'x-A^{-1}xz' = (A^2-A^{-2})y'$, 
$xx'-x'x=yy'-y'y=zz'-z'z=0$,\\ 
$Ax'y' -A^{-1}y'x' =(A^2-A^{-2})(z-A^{-1}(xy-x'y'))$,\\ 
$Ay'z'-A^{-1}z'y' = (A^2-A^{-2})x'$, 
$Az'x'-A^{-1}x'z' = (A^2-A^{-2})y'$, \\ 
and the long relation:\\ 
$A^6x^2 + A^{-2}y^2 +A^2z^2 + A^2x'^2 + A^{2}y'^2 +A^6z'^2 + A^2a^2 + 
A^4axx' + ayy'+A^4azz' +Ax'y'z - Axyz -Ax'yz' -A^5xy'z' -A^3axyz' 
-A^2(A^2 +A^{-2})^2$ 
\end{theorem} 
\begin{remark}\label{IX.9.3}
For $A=-1$ algebras of Proposition 5.5(2) and Theorem 5.6(1) are isomorphic.
Also algebras of Theorem 6.7 and 6.8 are isomorphic. It is the case because
for $A=-1$ the Kauffman bracket skein algebra depends only on the
fundamental group (it is explored in \cite{P-S-1,P-S-2} where an algebra is
associated to any group).\\
Of considerable interest is also the $R$ algebra monomorphism 
$$P: {\cal S}_{2,\infty}(F_{1,0}\times I;R[A^{\pm 1}],A) 
\to {\cal S}_{2,\infty}(F_{0,4}\times I;R[A^{\pm 1}],A)/ J$$ where 
$J$ is the ideal generated by expressions 
$(a_1-(A+A^{-1}),a_2-(A+A^{-1}),a_3-(A+A^{-1}), a_4 + (A+A^{-1}))$, and $P$ 
is given by $P(A)=A^2$, $P(x)=x_1$, $P(y)=x_2$ and $P(z)=x_3.$\\
If $\sqrt {-1} \in R$ then there is also another, related, $R$ algebra 
monomorphism $$P': {\cal S}_{2,\infty}(F_{1,0}\times I;R[A^{\pm 1}],A)         
\to {\cal S}_{2,\infty}(F_{0,4}\times I;R[A^{\pm 1}],A)/ J'$$ where     
$J'$ is the ideal generated by expressions 
$(a_1,a_2,a_3, a_4 - \sqrt {-1}(A^2-A^{-2}))$, and $P'$  
is given by $P'(A)=A^2$, $P'(x)=x_1$, $P'(y')=x_2$ and $P'(z')=x_3.$\\ 
We have noticed these monomorphisms by considering the 2-fold branched covering
of a torus over a 2-sphere (with the branching set of four points); \cite{B-P}.
\end{remark}

\section{Kauffman skein module ${\cal S}_{3,\infty}(M)$}\label{IX.10}
The deformation of the (unoriented) 2-move leads to the Kauffman 
skein relation and the Kauffman skein module.
\begin{definition}
Let $M$ be an oriented 3-manifold,
and put $R=\mathbb{Z}[a^{\pm1},x^{\pm1}]$.
The {\it Kauffman skein module} $S_{3,\infty}(M)$ of $M$
is defined to be
the $R$ module spanned by
unoriented framed links in $M$ ($R{\cal L}^{fr}$) subject 
to the relations:\\ 
the skein relation: $L_+ + L_- -x(L_0 + L_{\infty}) = 0$ (Fig. 11.1),\\
and the framing relation: $L^{(1)} = a L$.
\end{definition}

The Kauffman skein module has been computed for $S^3$, the solid torus, 
the product of a surface and the interval, and the projective space.
\begin{theorem}
\begin{enumerate}
\item[(i)] (Kauffman)\ \  $S_{3,\infty}(S^3) = R$ and $L = F_L(a,x)T_1 = 
\frac{a+ a^{-1} -x}{x}F_L(a,x)\emptyset$.
\item[(ii)] (Hoste-Kidwell-Turaev)\ \ $S_{3,\infty}(S^1\times D^2)$ is 
an infinitely generated free $R$-module.
\item[(iii)] (Lieberum)\ \ $S_{3,\infty}(F\times I)$ is a free $R$-module 
infinitely generated for $F$ not simply connected.
\item[(iv)] (Mroczkowski)\ \  $S_{3,\infty}(RP^3)$ is an infinitely generated 
free module with basis composed of standard unoriented
unlinks $L_n$ (links of Fig. 4.2 with orientation ignored).
\end{enumerate}
\end{theorem}

The relative Kauffman skein module has been found for a tangle, and
an annular tangle in the context of type $A$ and $B$ Brauer algebra 
\cite{Br}.
\begin{theorem}
\begin{enumerate}
\item[(i)] (Birman-Wenzl,J.Murakami,Morton-Traczyk)\ \ 
$S_{3,\infty}(D^2\times I;2n)$ is a free 
$(2n-1)(2n-3)\cdot...\cdot 3\cdot 1$- dimensional $R$-algebra.
\item[(ii)] (Goodman-Hauschild) $S_{3,\infty}(F_{0,2}\times I;2n)$ 
is a free $R$-algebra.
\end{enumerate}
\end{theorem}

\section{Skein module deformation of $3$- and $(2,2)$-moves}\label{IX.11}
In Sections 3-10 we
 gave  a description of skein
modules studied extensively until now. 
Their skein relations can be interpreted as a deformation of 
1- and 2-moves. More generally we can consider  
skein modules based on relations deforming n-moves:\ 
${\cal S}_n(M)= R{\cal L}/(b_0L_0 + b_1L_1 + b_2L_2 +...+b_{n-1}L_{n-1})$.
In the unoriented case, we can add to the relation the
term $b_{\infty}L_{\infty}$,
to get ${\cal S}_{n,\infty}(M)$, and also,
possibly, a framing relation. We can also consider skein modules 
based on skein relations deforming more general tangle moves e.g. 
rational moves. 

I will now describe two such examples which only recently
have been considered in more detail. The first example is based on
a deformation of the 3-move and the second on the deformation of
the $(2,2)$-move. The first one has been studied with
my students Tsukamoto and Veve. I denote the skein module
described in this example by
${\cal S}_{4,\infty}$ since it involves (in the skein relation)
4 horizontal positions and the vertical ($\infty$) smoothing.
\begin{definition}
Let $M$ be an oriented 3-manifold and let ${\cal L}_{fr}$ be
the set of unoriented
framed links in $M$ (including the empty link, $\emptyset$), and let
$R$ be any commutative ring with identity.  Then we define the
$(4,\infty)$ skein module as:
${\cal S}_{4,\infty}(M;R) = R{\cal L}_{fr}/I_{(4,\infty)}$,
where $I_{(4,\infty})$ is the submodule of $R{\cal L}_{fr}$
generated by the skein relation:\\
$b_0L_0  + b_1L_1 + b_2L_2 + b_3L_3 + b_{\infty}L_{\infty} = 0$
and
the framing relation:\\
$L^{(1)} = a L$ where $a,b_0,b_3$ are invertible elements in $R$ and
$b_1,b_2,b_{\infty}$ are any fixed elements of $R$ (see Fig. 11.1).
\end{definition}
\centerline{\psfig{figure=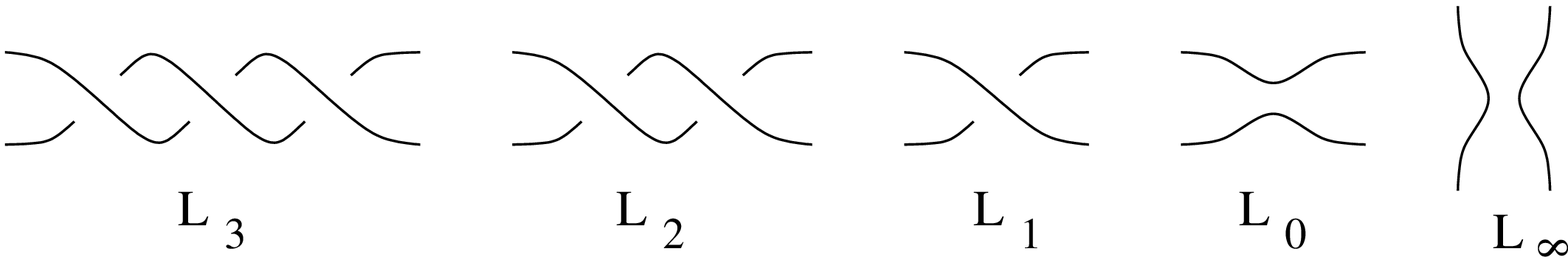,height=2.2cm}}
\centerline{Fig. 11.1.}
\ \\

It was conjectured for a while, generalizing the Montesinos-Nakanishi
3-move conjecture (Chapter I), that for $S^3$ the fourth skein 
module is generated by trivial links. However 
the counterexamples to the Montesinos-Nakanishi
3-move conjecture, Chapter VI, can be used to show that trivial links
``generically" do not generate ${\cal S}_{4,\infty}(S^3,R)$. 
\begin{proposition}\label{11.2} Assume that 
there is a proper ideal ${\cal I}\in R$ such that
$b_1,b_2 $ and $b_{\infty}$ are in ${\cal I}$. Then 
${\cal S}_{4,\infty}(S^3,R)$ is not generated by trivial links.
\end{proposition}
\begin{proof}
Let $R_I$ be the quotient ring $R/{\cal I}$. In $R_I$, $b_1=b_2=b_{\infty}=0$,
and $b_0$, $b_3$ and $a$ are invertible. Thus the skein relation reduces 
to $L_3 = (a_0a_3^{-1})L_0$. It follows from this observation that 
the family of links generates ${\cal S}_{4,\infty}(S^3,R_I)$ if 
and only if this family has a representative for every 3-move equivalence 
class. Thus ${\cal S}_{4,\infty}(S^3,R_I)$ is not generated by trivial links
and finally ${\cal S}_{4,\infty}(S^3,R)$ is not generated by trivial links.
The last conclusion follows immediately from our definitions but one 
can also observe that universal coefficient theorem can be used:
$${\cal S}_{4,\infty}(S^3,R)\otimes_RR/{\cal I}= 
{\cal S}_{4,\infty}(S^3,R_I).$$
\end{proof}

Similarly, our fourth skein module of $n$-tangles needs, generically,
 strictly more than $\prod_{i=1}^{n-1}(3^i+1)$ basic tangles (with
possible trivial components).  The lower bound was given in Section VI 
using interpretation of tangles as Lagrangian in the symplectic space 
of Fox 3-colorings.
In \cite{P-Ts} we analyzed extensively the possibility
that trivial links, $T_n$, are linearly independent.
This may happen if $b_{\infty} = 0$ and $b_0b_1=b_2b_3$.
These lead to the following conjecture:
\begin{conjecture}\label{11.3}
\begin{enumerate}
\item [\textup{(1)}] There is a polynomial invariant
of unoriented links in $S^3$,
$P_1(L) \in Z[x,t]$, which satisfies:
\begin{enumerate}
\item [\textup{(i)}]
Initial conditions: $P_1(T_n) = t^n$, where $T_n$ is a trivial
link of $n$ components.
\item [\textup{(ii)}]
Skein relation: $P_1(L_0) + xP_1(L_1) - xP_1(L_2) - P_1(L_3)=0$,
where $L_0,L_1,L_2,L_3$ is a standard, unoriented skein quadruple
($L_{i+1}$ is obtained from $L_{i}$ by a right-handed half-twist on
two arcs involved in $L_{i}$; compare Fig.3.3).
\end{enumerate}
\item [\textup{(2)}]
There is a polynomial invariant of unoriented framed links,
$P_2(L) \in Z[A^{\pm 1},t]$ which satisfies:
\begin{enumerate}
\item [\textup{(i)}] Initial conditions: $P_2(T_n) = t^n$,
\item [\textup{(ii)}]
Framing relation: $P_2(L^{(1)})=-A^3P_2(L)$ where $L^{(1)}$ is
obtained from a framed link $L$ by a positive half twist on its framing.
\item [\textup{(iii)}]
Skein relation: $P_2(L_0) + A(A^2 + A^{-2})P_2(L_1) +
(A^2 + A^{-2})P_2(L_2) + AP_2(L_3)=0$.
\end{enumerate}
\end{enumerate}
\end{conjecture}
The above conjectures assume that $b_{\infty}=0$ in our
skein relation. Let us consider, for a moment, the possibility that
$b_{\infty}$ is invertible in $R$. Using the ``denominator"
of our skein relation (Fig.11.2) we obtain the relation which
allows us to compute the effect of adding a trivial component
to a link $L$ (we write $t^n$ for the trivial link $T_n$):
$$(*)\ \ \ (a^{-3}b_3 + a^{-2}b_2 + a^{-1}b_1 + b_0 + b_{\infty}t)L=0$$
When considering the ``numerator" of the relation and its
mirror image (Fig. 11.2)
we obtain formulas for Hopf link summands, and because the unoriented
Hopf link is amphicheiral we can eliminate it from
our equations to get the following formula (**):
$$b_3(L\#H) + (ab_2 + b_1t +a^{-1}b_0 + ab_{\infty})L =0.$$
$$b_0(L\#H) + (a^{-1}b_1 +b_2t + ab_3 + a^2b_{\infty})L =0.$$
$$(**)\ \ \ ((b_0b_1 - b_2b_3)t + (a^{-1}b_0^2 -a b_3^2) +
(ab_0b_2 - a^{-1}b_1b_3) + b_{\infty}(ab_0 -a^2b_3))L = 0.$$
\ \\
\centerline{\psfig{figure=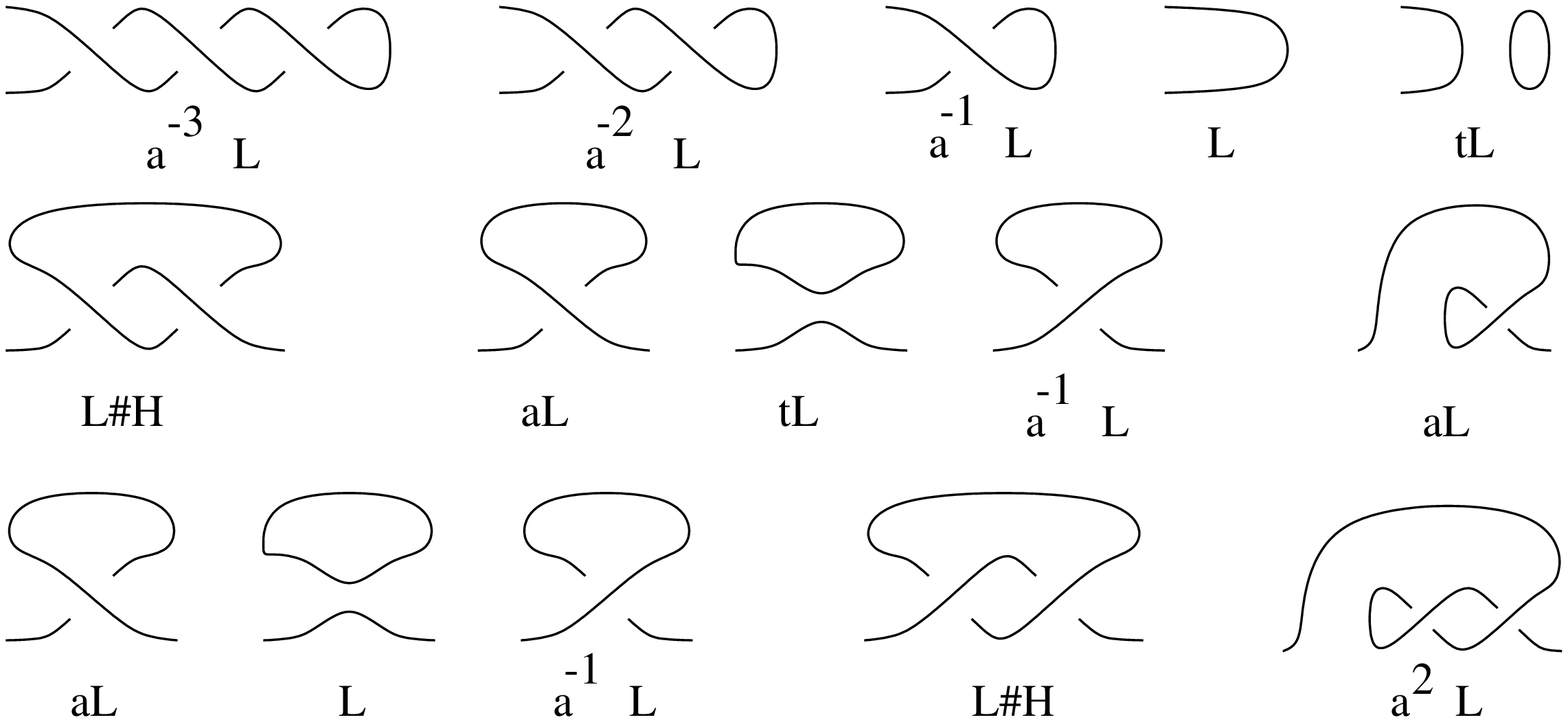,height=5.8cm}}
\centerline{Fig. 11.2}
\ \\
It is possible that ($*$) and ($**$) are the only relations in the module.
More precisely, we ask whether the submodule\footnote{We take into account
the fact that the Montesinos-Nakanishi
3-move conjecture does not hold \cite{D-P-1}.}
 of ${\cal S}_{4,\infty}(S^3;R)$
generated by trivial link
is the quotient ring $R[t]/({\cal I})$ where $t^i$ represents
the trivial link of $i$ components and ${\cal I}$ is the
ideal generated by ($*$) and ($**$) for $L=t$.\
The interesting substitution which satisfies the relations is
$b_0=b_3=a=1$, $b_1=b_2=x$, $b_{\infty}=y$. This may lead to
a new polynomial invariant (in $\mathbb{Z}[x,y]$) of unoriented links
in $S^3$ satisfying the skein relation
$L_3 +xL_2 + xL_1 + L_0 + yL_{\infty} = 0$.

What about the relations to the Fox colorings?
One such a relation, that was already mentioned,
is the use of 3-colorings to estimate the number of
basic n-tangles (by $\prod_{i=1}^{n-1} (3^i +1)$) for
the skein module ${\cal S}_{4,\infty}$. I am also convinced that
${\cal S}_{4,\infty}(S^3;R)$ contains full information about the space of
Fox 7-colorings. It would be a generalization of the fact that
the Kauffman bracket polynomial contains information about 3-colorings
and the Kauffman polynomial contains information about
5-colorings. In fact, Fran{\c c}ois Jaeger told me that he
knew how to form a short  skein relation
(of the type  $(\frac{p+1}{2},\infty)$) involving
 spaces of $p$-colorings.
Unfortunately,
Fran{\c c}ois died prematurely in 1997 and I do not know
how to prove his statement\footnote{If
$col_p(L)=|Col_p(L)|$ denotes the order of the space of Fox
$p$-colorings  of the link $L$, then
among $p+1$ links $L_0,L_1,...,L_{p-1}$, and $L_{\infty}$,
$p$ of them has the same order  $col_p(L)$
and one has its order $p$ times larger
\cite{P-20}. This leads to the relation of type $(p,\infty)$.
The relation between Jones polynomial (or the Kauffman bracket)
and $col_3(L)$ has the form: $col_3(L)= 3|V(e^{\pi i/3})|^2$ and
the formula relating the Kauffman polynomial and $col_5(L)$
has the form: $col_5(L)= 5|F(1,e^{2\pi i/5} + e^{-2\pi i/5})|^2$.
This seems to suggest that the formula discovered by Jaeger
involved Gaussian sums.}.

Finally, let us quickly  describe the skein module
related to the deformation of $(2,2)$-moves.
Because a $(2,2)$-move is equivalent to the rational
$\frac{5}{2}$-move, I will denote the skein module
by ${\cal S}_{\frac{5}{2}}(M;R)$.
\begin{definition}\label{11.4}
Let $M$ be an oriented 3-manifold. Let  ${\cal L}_{fr}$ be
 the set of unoriented framed links in $M$ (including the
empty link, $\emptyset$) and let $R$ be any commutative
ring with identity.  We define the $\frac{5}{2}$-skein module as
 ${\cal S}_{\frac{5}{2}}(M;R) = R{\cal L}_{fr}/(I_{\frac{5}{2}})$
where $I_{\frac{5}{2}}$ is the submodule of $R{\cal L}_{fr}$
generated by the skein relation:\\
\textup{(i)} \ \ $b_2L_2 + b_1L_1 + b_0L_0  + b_{\infty}L_{\infty} +
b_{-1}L_{-1} + b_{-\frac{1}{2}}L_{-\frac{1}{2}} = 0$,\\
 its mirror image:\\
$(\bar{i})$ \ \ $b'_2L_2 + b'_1L_1 + b'_0L_0  + b'_{\infty}L_{\infty} +
b'_{-1}L_{-1} + b'_{-\frac{1}{2}}L_{-\frac{1}{2}} = 0$\\
 and the framing relation:\\
$L^{(1)} = a L$, where $a,b_2,b_2',b_{-\frac{1}{2}},
b'_{-\frac{1}{2}}$ are invertible elements in $R$ and
$b_1,b'_1,b_0,b'_0$, $b_{-1},b'_{-1},b_{\infty}$,
and $b'_{\infty}$ are any fixed elements of $R$.
The links $L_2,L_1,L_0,L_{\infty},L_{-1},$ $L_{\frac{1}{2}}$
and $L_{-\frac{1}{2}}$
are illustrated in Fig. 11.3.\footnote{Our notation is based on
Conway's notation for rational tangles. However, it differs
from it by a sign change. The reason is that the Conway
convention for a positive crossing is generally not
used in the setting of skein relations.}
\end{definition}
\centerline{\psfig{figure=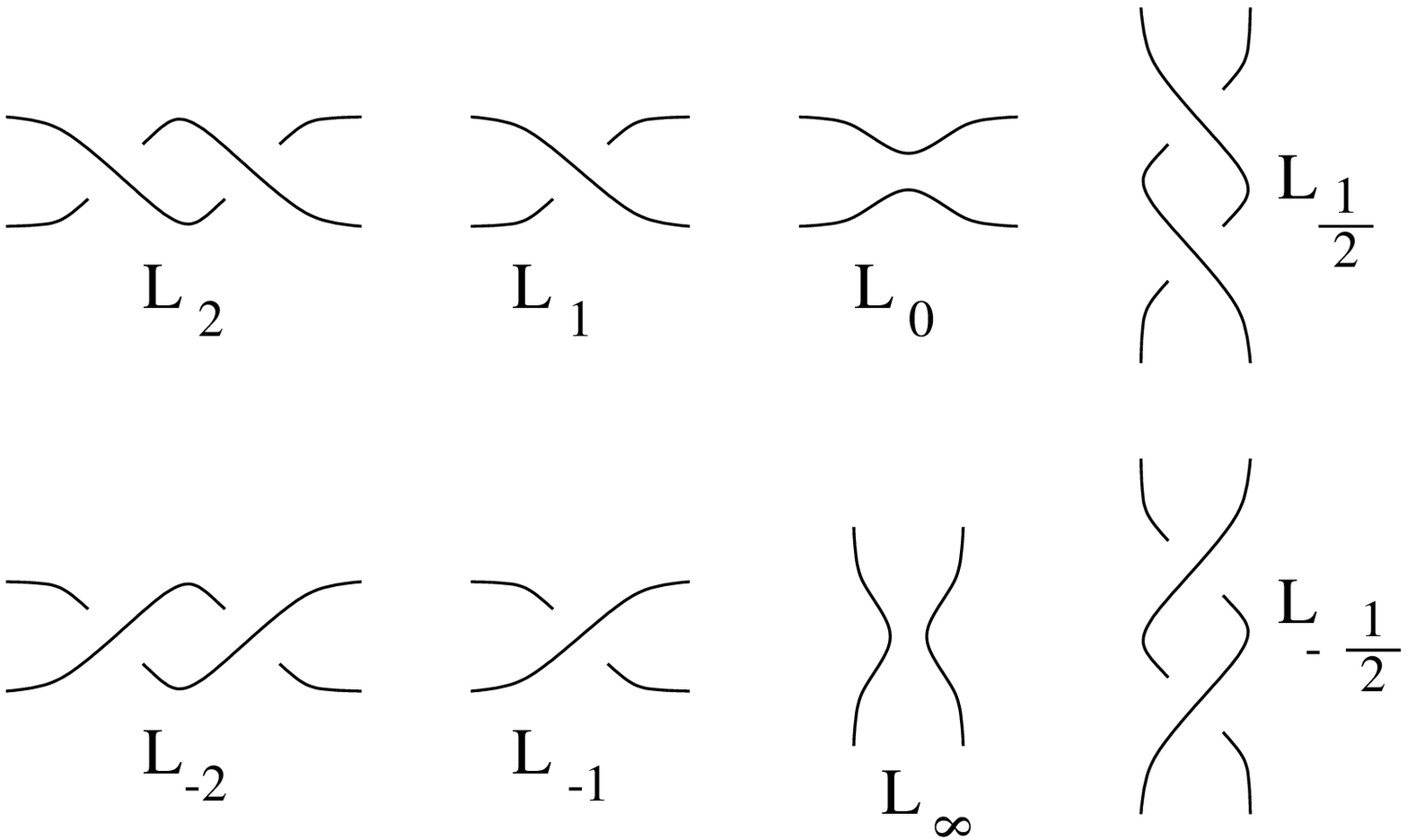,height=5.7cm}}
\centerline{Fig. 11.3}
\ \\

If we rotate the figure from the relation (i) we obtain:\\
(i')  $b_{-\frac{1}{2}}L_2 + b_{-1}L_1 + b_{\infty}L_0  + b_{0}L_{\infty} +
b_{1}L_{-1} + b_{2}L_{-\frac{1}{2}} = 0$\\
One can use (i) and (i') to eliminate $L_{-\frac{1}{2}}$
and to get the relation:\\
$(b_2^2- b^2_{-\frac{1}{2}}) L_2 +
(b_1b_2 - b_{-1}b_{-\frac{1}{2}})L_1 +
((b_0b_2 - b_{\infty}b_{-\frac{1}{2}})L_0 +
(b_{-1}b_2 - b_{1}b_{-\frac{1}{2}})L_{-1} +
(b_{\infty}b_2 - b_{0}b_{-\frac{1}{2}})L_{\infty}=0$.\\
Thus, either we deal with the shorter relation (essentially the one in
the fourth skein module described before) or all
coefficients are equal to 0 and therefore (assuming
that there are no zero divisors in $R$)
$b_2 = \varepsilon b_{-\frac{1}{2}}$, $b_1 = \varepsilon b_{-1}$,
and $b_0 = \varepsilon b_{\infty}$. Similarly, we would get:
$b'_2 = \varepsilon b'_{-\frac{1}{2}}$, $b'_1 = \varepsilon b'_{-1}$,
and $b'_0 = \varepsilon b'_{\infty}$, where $\varepsilon = \pm 1$.
Assume, for simplicity, that $\varepsilon =1$. Further relations
among coefficients follow from the computation of the
Hopf link component using the amphicheirality of the
unoriented Hopf link.
Namely, by comparing diagrams in Figure 3.6 and their mirror images we get
$$L\#H= -b_2^{-1}(b_1(a + a^{-1}) + a^{-2}b_2 + b_0(1+T_1))L$$
$$L\#H= -{b'}_2^{-1}(b'_1(a + a^{-1}) + a^{2}b'_2 + b'_0(1+T_1))L.$$
Possibly, the above equalities give the only other relations among
coefficients (in the case of $S^3$). I would present below
the simpler question (assuming $a=1, b_x=b'_x$ and writing
$t^n$ for $T_n$).\\
\ \\
\centerline{\psfig{figure=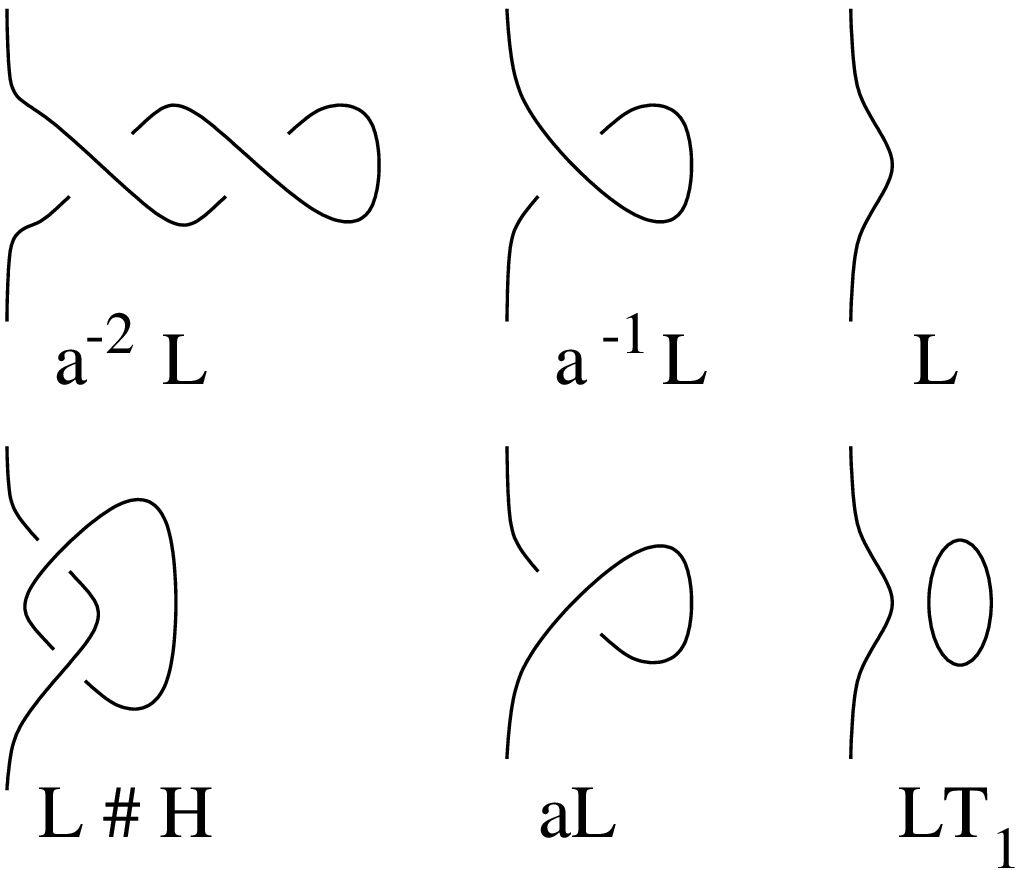,height=6.4cm}}
\centerline{Fig. 11.4}

\begin{question}\label{11.5}
Is there a polynomial invariant of unoriented links in $S^3$,
$P_{\frac{5}{2}}(L) \in {\mathbb Z}[b_0,b_1,t]$, which satisfies the
following conditions?
\begin{enumerate}
\item [\textup{(i)}]
Initial conditions: $P_{\frac{5}{2}}(T_n) = t^n$, where $T_n$ is a trivial
link of $n$ components.
\item [\textup{(ii)}]
Skein relations
$$P_{\frac{5}{2}}(L_2) + b_1P_{\frac{5}{2}}(L_1) +  b_0P_{\frac{5}{2}}(L_0) +
b_0P_{\frac{5}{2}}(L_{\infty}) + b_1P_{\frac{5}{2}}(L_{-1}) +
P_{\frac{5}{2}}(L_{-\frac{1}{2}})
=0.$$
$$P_{\frac{5}{2}}(L_{-2}) + b_1P_{\frac{5}{2}}(L_{-1}) +
b_0P_{\frac{5}{2}}(L_0) +
b_0P_{\frac{5}{2}}(L_{\infty}) + b_1P_{\frac{5}{2}}(L_{1}) +
P_{\frac{5}{2}}(L_{\frac{1}{2}})
=0.$$
\end{enumerate}
\end{question}
Notice that by taking the difference of our skein relations one gets
the interesting identity:
$$P_{\frac{5}{2}}(L_2)  - P_{\frac{5}{2}}(L_{-2}) =
P_{\frac{5}{2}}(L_{\frac{1}{2}}) - P_{\frac{5}{2}}(L_{-\frac{1}{2}}).$$
Nobody has yet studied the skein module
${\cal S}_{\frac{5}{2}}(M;R)$ seriously so everything that you can find will
be a new research, even a table of the polynomial $P_{\frac{5}{2}}(L)$
for small links, $L$.\\

\section{Vassiliev-Gusarov skein modules}\label{IX.12}

Let ${\cal K}^{sg}$ denote the set of singular oriented knots in $M$ where 
we allow only immersions of $S^1$ with, possibly, double points, up to 
ambient isotopy; additionally for any double point we choose orientation 
of a small ball around it (if $M$ is oriented the chosen orientation of 
the ball agrees with that of $M$). 
In $R{\cal K}^{sg}$ we 
consider resolving singularity relations $\sim$: $K_{cr} = K_+ - 
K_-$ ; see Fig. 12.1. 
\ \\
\centerline{\psfig{figure=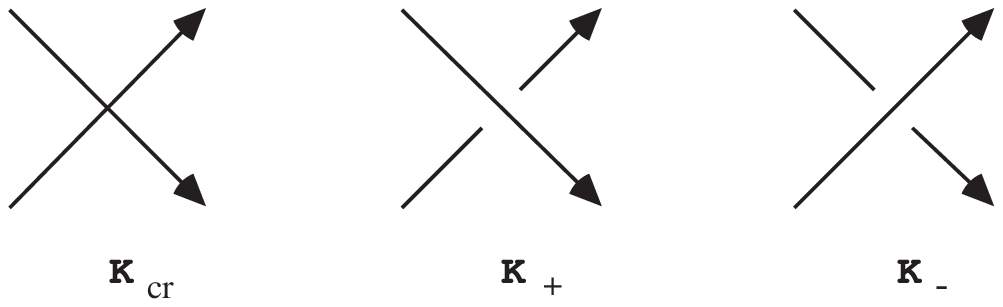}} 

\begin{center} 
Figure 12.1 
\end{center} 

$R{\cal K}^{sg}/\sim$ is obviously $R$-isomorphic 
to $R\cal K$. Let $C_m$ be a 
submodule of $R{\cal K}^{sg}/\sim =R\cal K$ generated by immersed knots 
with $m$ double points. 
The $m$'th Vassiliev-Gusarov 
skein module $W_m(M,R)$ is defined
by $W_m(M,R) = R{ \cal K} / C_{m+1}$. \ \ 
We have the filtration:\\ 
 $...\subset C_m \subset ...\subset C_1 \subset C_0 = R{\cal 
 K}$ 
and therefore we have also 
the sequence of epimorphisms  $\{1\} \leftarrow W_0 \leftarrow 
W_1 \leftarrow W_2 \leftarrow ...\leftarrow W_m 
\leftarrow ...$ 
We define the V-G skein module $W_{\infty}(M;R)$ as the inverse limit 
$W_{\infty}(M;R) = 
lim_{\leftarrow}W_m(M;R).$ 
Equivalently the V-G skein module is the completion $\widehat {R\cal K}$ 
of $R\cal K$ with respect to 
the topology yielded by the sequence of descending submodules $C_i$. 
The kernel of the natural $R$-homomorphism $\theta : R{\cal K} \to 
\widehat {R{\cal K}}$ is equal to $\bigcap_{i=0}^{\infty}C_i$. 

$W_{\infty}(S^3;R)$ is a Hopf algebra with: $\mu(K_1,K_2)=K_1\#K_2$,
$i(1)=T_1$, $\Delta(K)=K\otimes K$, $\epsilon (K)=1$,\ 
$S(K)=T_1 - (K-T_1) + (K-T_1)^2 - (K-T_1)^3 +...$, where $K_1,K_2$ and 
$K$ are non-singular knots; \cite{P-9}. 

A V-G invariant of degree $m$ of knots is defined as an element of the 
dual space $V^m(M,R)= W^*_m(M,R) = Hom(W_m(M,R),R)$ 
(sometimes it is defined more generally,
as an element of $Hom_Z(W_m(M,Z),A)$, where $A$ is an abelian group).\\

We can modify the definition of V-G skein modules by resolving 
the singular points differently, as suggested in \cite{P-9}. 
The most natural resolutions are those 
suggested by Jones type skein relations. For example the resolution
$L_{cr} -v^{-1}L_+ -v L_- -zL_0$ proposed in \cite{P-9} was analyzed 
in more detail by Y. Rong and R. Lickorish \cite{Ron-2,L-R}. 
The resolution of a crossing of an 
unoriented framed singular link $L_{cr} - L_+ -AL_0 -A^{-1}L_{\infty}$ 
was hinted in \cite{B-F-K-2}. Further work in this direction was 
done in \cite{A-T-1,A-T-2}.

\section{A glimpse into 4-dimensional skein modules}\label{IX.13}
The following approach to 4-dimensional skein modules is based on ideas
of Kamada and Kawauchi \cite{Kaw-3}.
\begin{definition}
Let $M$ be an oriented 4-dimensional manifold and ${\cal F}$ denote the set
of immersed surfaces in $M$ with possible double points, up to ambient
isotopy. Let $R$ be a commutative ring with unit 
and let $R\cal F$ denote the free $R$-module generated by $\cal F$.
Consider the submodule $K$ of $R\cal F$ generated by expressions:
$r_0F_0 + r_1F_1 +r_2F_2$ where $r_i \in R$ and $F_1$ and $F_2$ are
obtained from $F_0 \in {\cal F}$ as follows (Fig. 13.1):\\
Let $\gamma$ be any curve joining two points on $F_0$, otherwise 
disjoint from $F_0$. 
\begin{enumerate}
\item[(1)] $F_1$ is obtained from $F_0$ by performing index 1 surgery 
on $F_1$ along $\gamma$. 
\item[(2)] $F_2$ is obtained from $F_0$ by taking, instead of a 1-surgery, 
an immersed surface with two additional double points (one should carefully 
choose a convention - there are two choices).
\end{enumerate}
Now in the quotient module, ${\cal K}=R{\cal F}/K$, 
consider the filtration $\{C_i\}$, where
$C_i$ is generated by surfaces from ${\cal F}$ with $i$ double
points. The completion of the space with respect to the filtration
is called the KK skein module of a four dimensional manifold $M$ and denoted
by $KK_{\infty}(M;R)$. The quotient of ${\cal K}$  by $C_{i+1}$ are 
degree $i$, KK
skein modules, $KK_i(M;R)$. Dual elements (or, more generally
homomorphisms to any abelian group) are called degree $i$, KK invariants
of a 4-manifold.
\end{definition}

\ \\
\ \\
\centerline{\psfig{figure=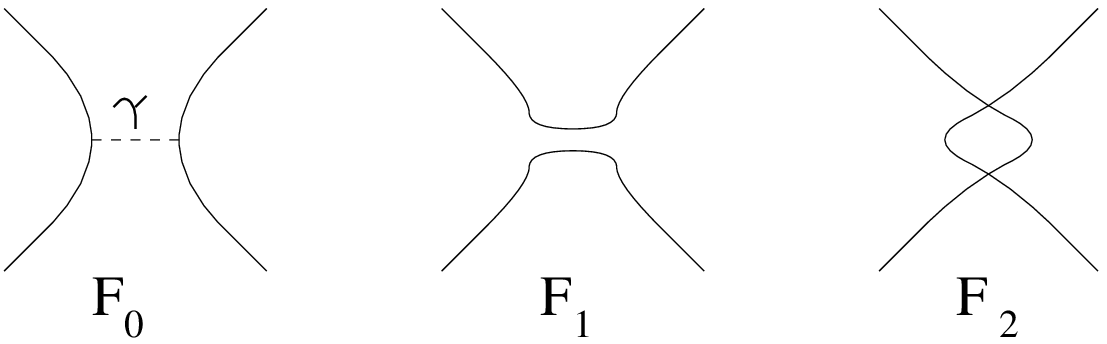}}
\begin{center}
Fig.~13.1
\end{center}

\ \\ \ \\ \ \\
\noindent \textsc{Dept. of Mathematics, Old Main Bldg., 1922 F St. NW \\
The George Washington University, Washington, DC 20052}\\
e-mail: {\tt przytyck@gwu.edu}


\begin{thebibliography}{99 1000-2000}

\bibitem [A-V]{A-V}
J.~W.~Alexander, Letter to Oswald Veblen, 1919,
Papers of Oswald Veblen, 1881-1960 (bulk 1920-1960),
Archival Manuscript Material (Collection), Library of Congress.

\bibitem[Al-3]{Al-3} 
J.~W.~Alexander, Topological invariants of knots and links,
{\it Trans. Amer. Math. Soc.} 30, 1928, 275-306.

\bibitem[A-T-1]{A-T-1}
J.~E. Andersen, V. Turaev, Higher skein modules, 
{\it Journal of Knot Theory and its Ramifications}, 8(8), 1999, 963--984.\\
e-print: http://front.math.ucdavis.edu/math.GT/9812071 

\bibitem[A-T-2]{A-T-2}
J.~E. Andersen, V. Turaev,
Higher skein modules II,  in ``Topology, Ergodic Theory, Real Algebraic 
Geometry. Rokhlin's Memorial" {\it Amer. Math. Soc. 
Transl.} ser 2, vol 202, 2001.

\bibitem[Ba-1]{Ba-1} D. Bar-Natan,
 On the Vassiliev knot invariants, {\em Topology}, 
34(2), 1995, 423-472.

\bibitem[Birk-1]{Birk-1} 
G.~D. Birkhoff,  A determinant formula for the number of ways 
      of coloring a map, {\em Ann. of Math.} (2)14, 1912,  42-46.

\bibitem[Bi-2]{Bi-2}
J.~S. Birman, On the Jones polynomial of closed 3-braids, {\it Invent. 
Math.}, 81(2), 1985, 287-294. 

\bibitem[B-W]{B-W}
J.~S. Birman, H. Wenzl, Braids, link polynomials and a new algebra, 
{\it Trans. Amer. Math. Soc.}, 313(1), 1989, 249-273.

\bibitem [Bl] {Bl} 
T.~S. Blyth, {\it Module theory}, Clarendon Press, Oxford 1977. 

\bibitem [Br]{Br}
R. Brauer, On algebras which are connected with semisimple Lie groups, 
{\it Ann. of Math.} 38, 1937, 857-872.

\bibitem[Bu-1]{Bu-1}
 D. Bullock, Skein related links in 3-manifolds, {\it Topology and
its applications}, 60(3), 1994, 235-248.

\bibitem[Bu-2]{Bu-2} D. Bullock, The ($2,\infty$)--skein module of the
complement of a ($2,2p+1$) torus knot, {\it Journal of
Knot theory}, 4(4) 1995, 619-632.

\bibitem[Bu-3]{Bu-3} 
D. Bullock, A finite set of generators for the Kauffman bracket
skein algebra, {\it Math. Z.} 231(1), 1999, 91--101
(in the preprint form (1995) the paper had the title:
An integral invariant of 3-manifolds derived from the
Kauffman bracket).

\bibitem[Bu-4]{Bu-4} D. Bullock, 
Estimating a Skein Module with $SL_2(C)$ characters, 
{\it Proc.  Amer. Math. Soc.} 125(6), 1997, 1835--1839.

\bibitem[Bu-5]{Bu-5} D. Bullock, 
Estimating the States of the Kauffman Bracket Skein Module, 
{\it Knot theory} (Warsaw, 1995), 23--28, Banach Center Publ., 42,
Polish Acad. Sci., Warsaw, 1998.

\bibitem[Bu-6]{Bu-6}
D. Bullock, Rings of $SL_2(C)$ characters and the Kauffman
bracket skein module, {\it Comment. Math. Helv.} 72(4),
1997, 521--542.

\bibitem[Bu-7]{Bu-7}
D. Bullock, On the Kauffman bracket skein module of surgery on 
a trefoil, {\it Pacific J. Math.}, 178(1), 1997, 37--51.

\bibitem[B-F-K-1]{B-F-K-1}
D. Bullock, C. Frohman, J. Kania-Bartoszy\'nska,
Skein homology, {\it Canad. Math. Bull.} 41(2), 1998, 140--144.\\
e-print:\ http://arxiv.org/abs/q-alg/9701019

\bibitem [B-L] {B-L} D. Bullock, W. Lofaro,
The Kauffman bracket skein module of a twist knot exterior.\\
e-print:\ 
http://front.math.ucdavis.edu/math.QA/0402102

\bibitem[B-F-K-2]{B-F-K-2}
D. Bullock, C. Frohman, J. Kania-Bartoszy\'nska, 
Understanding the Kauffman bracket skein module,
{\it J. Knot Theory Ramifications} 8(3), 1999, 265--277.\\
e-print: http://arxiv.org/abs/q-alg/9604013

\bibitem[B-P]{B-P} 
D. Bullock,  J.~H.~Przytycki, Multiplicative structure of Kauffman bracket 
skein module quantizations, 
{\it Proc. Amer. Math. Soc.}, 128(3), 2000, 923--931
(in the preprint form the paper had the title:
Kauffman bracket skein module quantization 
of symmetric algebra and $so(3)$).\\
e-print:\ http://arxiv.org/abs/math.QA/9902117

\bibitem  
[C-E] {C-E}  
H. Cartan, S. Eilenberg, {\it Homological Algebra}, 
Princeton University Press,  1956.  

\bibitem[Ca]{Ca}  A. Caudron, Classification des noeuds et des enlacements,
             Prepublications, Univ. Paris-Sud, Orsay 1981.

\bibitem[Cher]{Cher}
V. Chernov, Framed knots in 3-manifolds and affine self-linking numbers;\\
e-print: http://front.math.ucdavis.edu/math.GT/0105139

\bibitem[Co-1]{Co-1} 
J.~H.~Conway, An enumeration of knots and links, 
{\em Computational problems in abstract algebra} (ed. J.Leech), 
Pergamon Press, 1969, 329 - 358.

\bibitem[Co-2]{Co-2} 
J.~H. Conway, Lecture, University of Illinois at Chicago, spring 1978. 

\bibitem[Co-3]{Co-3} 
J.~H. Conway, Talks at Cambridge Math. Conf., summer, 1979.

\bibitem[Cox]{Cox} 
H.~S.~M. Coxeter, Factor groups of the braid group, Proc. Fourth Canadian  
Math. Congress, Banff, 1957, 95-122.  

\bibitem [D-I-P]{D-I-P}
M.~K.~D{\c a}bkowski, M.~Ishiwata, J.~H. Przytycki,
$(2,2)$-move equivalence for 3-braids and 3-bridge links,
in preparation.

\bibitem [D-P-1]{D-P-1}
M.~K. D{\c a}bkowski, J.~H. Przytycki, Burnside obstructions
to the Montesinos-Nakanishi 3-move conjecture,
{\it Geometry and Topology}, June, 2002, 335-360.\\
e-print:\  http://front.math.ucdavis.edu/math.GT/0205040

\bibitem [D-P-2]{D-P-2}
M.~K. D{\c a}bkowski, J.~H. Przytycki,
Unexpected connection between Burnside groups and Knot Theory, 
{\it PNAS}, to appear.  \\ 
e-print:\ http://front.math.ucdavis.edu/math.GT/0309140

\bibitem [D-P-3]{D-P-3}
M.~K. D{\c a}bkowski, J.~H. Przytycki,
Burnside groups in knot theory, preprint 2004.

\bibitem[D-W]{D-W}
R. Dedekind, H. Weber, Theorie der algebraischen Funktionen 
einer Ver\"anderlichen, {\em J. de Crelle}, XCII, 1882, 181-290. 
(Also in R.Dedekind, Ges.Math.Werke,t. I, p. 248-349). 

\bibitem[Die]{Die} 
J.~Dieudonn\'e, {\em A history of algebraic and differential topology
1900-1960}, Birkh\"auser, Boston, Basel, 1989.

\bibitem [Dr]{Dr}
J.~V. Drobotukhina, An analogue of the Jones polynomial for links 
in ${\R}P^3$ and a generalization of the kauffman-Murasugi theorem,
{\it Algebra i Analiz} 2(3), 1990, 171-191; English transl., 
{\it Leningrad Math. J.}, 2(3), 1991, 613-630.

\bibitem[FYHLMO] {FYHLMO} 
P.~Freyd, D.~Yetter, J.~Hoste, W.~B.~R.~Lickorish, K.~Millett,
A.~Ocneanu, A new polynomial invariant of knots
and links, {\em Bull. Amer. Math. Soc.}, 12, 1985, 239-249.

\bibitem[G-J]{G-J}
M.~R.~Garey, D.~S. Johnson, {\em Computers and Intractability: 
a Guide to Theory of NP Completeness}, Ed. V.Klee, New York, 1979.

\bibitem[Gi]{Gi}  C.~A. Giller, A Family of links and the Conway calculus,
{\it Trans. Amer. Math. Soc.}, 270(1), 1982, 75-109.

\bibitem[G-K-P]{G-K-P}
P. Gilmer, J. Kania-Bartoszy\'nska, J.~H. Przytycki,
3-manifold invariants and periodicity of rational homology spheres,
{\it Algebr. Geom. Topol.}, 2, 2002, 825-842.\\
http://arxiv.org/abs/math.GT/9807011

\bibitem[Gol] {Gol}
W. Goldman, Invariant functions of Lie groups and Hamiltonian flows
of surface group representation, {\it Invent. Math.}, 85, 1986, 263-302.

\bibitem[G-H]{G-H}
F.~M. Goodman, H.~M. Hauschild, Affine Birman-Wenzl-Murakami Algebras 
and Tangles in the Solid Torus, {\it Fundamenta Mathematicae}, 
to appear.\\

\bibitem[Ho-K]{Ho-K} J. Hoste, M. Kidwell, Dichromatic link invariants, 
{\em Trans.  Amer. Math. Soc.}, 321(1), 1990, 197-229; see also
the preliminary version of this paper: ``Invariants of colored links",
preprint, March 1987.

\bibitem[H-P-1]{H-P-1} J.~Hoste, J.~H. Przytycki, 
An invariant of dichromatic links,
{\it Proc. Amer. Math.Soc.}, 105(4) 1989, 1003-1007.

\bibitem[H-P-2]{H-P-2} 
J.~Hoste, J.~H. Przytycki, Homotopy skein modules of orientable
3-manifolds, {\it Math. Proc. Camb. Phil. Soc.}, 108, 1990, 475-488.

\bibitem[H-P-3]{H-P-3} 
J.~Hoste, J.~H. Przytycki, A survey of skein modules of 3-manifolds,
in  Knots 90, Proceedings of the International Conference on Knot
Theory and Related Topics, Osaka (Japan), August 15-19, 1990, Editor
A.~Kawauchi, Walter de Gruyter 1992, 363-379.

\bibitem[H-P-4]{H-P-4} 
J.~Hoste, J.~H. Przytycki, The $(2,\infty)$-skein module of lens
 spaces; a generalization of the Jones polynomial,
{\it Journal of Knot Theory and Its Ramifications}, 2(3), 1993, 321-333.

\bibitem[H-P-5]{H-P-5}
J.~Hoste, J.~H. Przytycki, The skein module of genus 1 Whitehead
type manifolds, {\em Journal of Knot Theory and Its Ramifications},
4(3), 1995, 411-427.

\bibitem[H-P-6]{H-P-6}
J.~Hoste, J.~H. Przytycki, The Kauffman bracket skein module of
$S^1 \times S^2$, {\em Math. Z.},  220(1), 1995, 63-73.

\bibitem[Hud]{Hud}
J.~F.~P.~Hudson, {\it Piecewise linear topology}, Benjamin Inc. N.Y.,1969.

\bibitem[Iw]{Iw}
N. Iwahori, On the structure of a Hecke ring of a Chevalley group over
a finite field, {\it J.~Fac.~Sci.~Univ. Tokyo Sect.~I}, 10, 1964, 215-236.

\bibitem[Ja-1]{Ja-1}
F. Jaeger, On Tutte polynomials and link polynomials, {\em Proc. Amer. 
Math. Soc.}, 103(2), 1988, 647-654. 
 
\bibitem[Ja-2]{Ja-2}  
F. Jaeger, Composition products and models for the Homfly
polynomial, {\em L'Enseignement Math\'ematique}, 35, 1989, 323-361.

\bibitem[J-V-W]{J-V-W} 
 F. Jaeger, D.~L. Vertigan, D.~J.~A. Welsh, On the Computational  
Complexity of the Jones and Tutte Polynomials,  
{\em Math. Proc. Camb. Phil. Soc.}, 108, 1990, 35-53. 

\bibitem[Jo-3]{Jo-3}
 V.~F.~R. Jones. On a certain value of the Kauffman polynomial.
{\em Comm. Math. Phys.}, 125(3), 1989, 459--467.

\bibitem[Kai-1]{Kai-1}
U. Kaiser, Link homotopy and skein modules of 3-manifolds,
in {\em Geometric Topology}, ed.: C.Gordon, Y.Moriah, B.Wajnryb,
Contemporary Mathematics 164, 1994, 59-77.

\bibitem[Kai-2]{Kai-2}
U. Kaiser, Presentations of homotopy skein modules of oriented
3-manifolds, {\it J. Knot Theory Ramifications} 10(3), 2001, 461--491.\\
e-print:\ http://front.math.ucdavis.edu/math.GT/0007031

\bibitem[Kai-3]{Kai-3}
U. Kaiser, 
Quantum deformations of fundamental groups of oriented 3-manifolds,
{\it Trans. Amer. Math. Soc.}, November 25, 2003.\\
e-print;\ http://front.math.ucdavis.edu/math.GT/0202293 

\bibitem[Kai-4]{Kai-4}
U. Kaiser, Link theory in manifolds, Lecture Notes in Mathematics, i
Vol. 1669, 1997,  Springer Verlag, pp. 169. 

\bibitem[Ka-1]{Ka-1}  
L.~H. Kauffman, Combinatorics and knot theory, Contemporary
Math., Vol.20, 1983, 181-200.

\bibitem[Ka-3]{Ka-3} 
L.~H. Kauffman, {\em On knots}, Annals of Math. Studies, 115,
Princeton University Press, 1987.

\bibitem[Ka-8]{Ka-8}
L.~H. Kauffman, {\it Formal knot theory}, Mathematical Notes 30,
Princeton University Press, 1983.

\bibitem[Kaw-3]{Kaw-3}
A.~Kawauchi, Personal communication, July 1996.

\bibitem[Kir]{Kir}
 R. Kirby, Problems in low-dimensional topology; Geometric Topology  
(Proceedings of the Georgia International Topology Conference, 1993),  
Studies in Advanced Mathematics, Volume 2 part 2., Ed. W.Kazez, AMS/IP,  
1997, 35-473. 

\bibitem[L-P]{L-P}
S. Lambropoulou, J.~H.~Przytycki,
Hecke algebra approach to skein modules of lens spaces (with S.Lambropoulou);
in preparation.

\bibitem [Le]{Le}
T. Le, The colored Jones polynomial and the A-polynomial of 2-bridge knots,
presentation at the {\it Mini-Conference in Logic and Topology},
GWU, April 9, 2004.

\bibitem[Li-10]{Li-10} 
W.~B.~R.~Lickorish, The panorama of polynomials for knots,
links and skeins. In {\em Braids}, ed. J.S.Birman and A.L.Libgober,
Contemporary Math. Vol. 78, 1988, 399-414.

\bibitem[Li-3]{Li-3} 
W.~B.~R. Lickorish.  Polynomials for links.{\em  Bull. London Math.
 Soc.}, 20, 1988, 558-588.

\bibitem[L-M-1]{L-M-1}  
W.~B.~R.~Lickorish, K.~Millett, A polynomial invariant
of oriented links, {\it Topology}, 26, 1987, 107-141.

\bibitem[L-M-3]{L-M-3} 
W.~B.~R.~Lickorish, K.~C.~Millett, Some evaluations of link
polynomials, {\it Comment. Math. Helv.}, 61(1986), 349-359.

\bibitem[L-M-4]{L-M-4}
 W.~B.~R.~Lickorish, K.~C.~Millett, An evaluation of the F-polynomial
of a link, {\em Differential topology}, Proc. 2nd Topology
Symp., Siegen / FRG 1987, Lect. Notes Math. 1350, 1988, 104-108.

\bibitem[L-R]{L-R}
 W.~B.~R. Lickorish, Y. Rong, 
On Derivatives of Link Polynomials,
{\it Topology and its Applications}, 87(1), 1998, .

\bibitem[Lieb]{Lieb}
J. Lieberum, Skein modules of links in cylinders over surfaces,
{\it International Journal of mathematics and Math. Sci.}, 32(9), 2002,
 515-554.\\ 
e-print:\ arXiv:math.QA/9911174

\bibitem[Mas]{Mas}
G. Masbaum, The spin refined Kauffman bracket skein module
of $S^1\times S^2$ and lend spaces,
{\it Manuscripta Math.} 91(4), 1996, 495--509.

\bibitem[Mc]{Mc}
D. McCullough, Mappings of reducible 3-manifolds; In 
{\it Geometric and Algebraic Topology}, pp. 61-76, Banach Center 
Publications, Vol. 18, PWN, 1986.

\bibitem [Mo-Tr-2]{Mo-Tr-2}
H.~R. Morton, P. Traczyk, Knots and algebras, {\it Contribuciones Matematicas
en homenaje al profesor D.Antonio Plans Sanz de Bremond},
ed. E.Martin-Peinador and A.Rodez Usan, University of Zaragoza,
1990, 201-220.

\bibitem[Mro]{Mro}
M.~ Mroczkowski, Polynomial Invariants of links in the projective space,
{\it Fundamenta Mathematicae}, 184, Proceedings of Knots in Poland 2003, 
December 2004, to appear.\\
http://front.math.ucdavis.edu/math.GT/0312205


\bibitem[Mur-1]{Mur-1} 
H. Murakami, A recursive calculation of the Arf invariant of a link,
{\em J. Math. Soc. Japan.}, 38, 1986, 335-338.

\bibitem [Mura-2]{Mura-2}
J. Murakami, The Kauffman polynomial of links and representation theory,
{\it Osaka J. Math.}, 24, 1987, 745-758.

\bibitem[Po]{Po} H.~Poincar\'e,   
Analysis Situs (\&12), {\em  Journal d'Ecole Polytechnique 
Normale}, 1, 1895,  1-121. 

\bibitem [P-1]{P-1}
J.~H. Przytycki, Survey on recent invariants  in classical
knot theory, Uniwersytet Warszawski, Preprinty 6,8,9; 1986.

\bibitem[P-5]{P-5} J.~H. Przytycki, Skein modules of 3-manifolds,
{\em Bull. Ac. Pol.: Math.}; 39(1-2), 1991, 91-100.

\bibitem[P-6]{P-6}
J.~H. Przytycki,
Homotopy and q-homotopy skein modules of 3-manifolds: an example in
Algebra Situs; In:
{\it Knots, Braids, and Mapping Class Groups: Papers
dedicated to Professor Joan Birman},
Ed. J. Gilman, W. Menasco, and X.-S. Lin, International Press.,
AMS/IP Series on Advanced Mathematics, Vol 24,
Co., Cambridge, MA, 2001, 143-170.
(extended version of: $q$-analog of the homotopy skein module,
preprint, Knoxville 1991).\\
e-print:\ http://arxiv.org/abs/math.GT/0402304.

\bibitem[P-9]{P-9} J.~H.~Przytycki,
Vassiliev-Gusarov skein modules of 3-manifolds and criteria for periodicity
of knots, Low-Dimensional Topology, Knoxville, 1992 ed.: Klaus Johannson,
International Press Co., Cambridge, MA 02238, 1994, 157-176.

\bibitem[P-12]{P-12}
J.~H. Przytycki, A q-analogue of the first homology group of a 3-manifold,
{\it Contemporary Mathematics} 214, Perspectives on Quantization
(Proceedings of the joint AMS-IMS-SIAM conference on Quantization,
Mount Holyoke College, 1996); Ed. L.A.Coburn, M.A.Rieffel, AMS 1998,
135-144.

\bibitem[P-14]{P-14}
J.~H.~Przytycki, Skein module of links in a handlebody, 
Topology 90, Proc. of the Research Semester 
in Low Dimensional Topology at OSU, Editors: B.Apanasov, 
W.D.Neumann, a.W.Reid, 
L.Siebenmann, De Gruyter Verlag,  1992; 315-342. 

\bibitem[P-15]{P-15} J.~H.~Przytycki,
Quantum group of links in a handlebody {\it Contemporary Math:
Deformation Theory and Quantum Groups with Applications
to Mathematical Physics}, M.Gerstenhaber  and J.D.Stasheff, Editors,
Volume 134, 1992, 235-245.

\bibitem[P-18]{P-18}
J.~H. Przytycki,
{\it Teoria w\c ez\l\'ow: podej\'scie kombinatoryczne},
(Knots: combinatorial approach to the knot theory),
Script, Warsaw, August 1995, 240+XLVIIIpp.

\bibitem[P-19]{P-19} 
J.~H. Przytycki, Algebraic topology based on knots: an introduction, 
{\it Knots 96}, Proceedings of the Fifth International Research Institute
of MSJ, edited by Shin'ichi Suzuki, World Scientific Publishing Co., 1997,
279-297.

\bibitem[P-20]{P-20}
J.~H.~Przytycki, 3-coloring and other elementary
invariants of knots, Banach Center Publications, Vol. 42,
{\it Knot Theory}, 1998, 275-295.

\bibitem[P-30]{P-30}
J.~H. Przytycki,
{\it Algebraic topology based on knots},
Series on Knots and Everything - Vol. 18, World Scientific, in preparation.

\bibitem[P-S-1]{P-S-1} J.~H.~Przytycki, A.~S.~Sikora, 
Skein algebra of a group, , Banach Center Publications,
Vol. 42, {\it Knot Theory}, 1998, 297-306.

\bibitem[P-S-2]{P-S-2} J.~H.~Przytycki, A.~S.~Sikora,
On skein algebras and $Sl_2(C)$-character varieties, 
{\it Topology}, 39(1), 2000, 115-148.\\
e-print: \ http://arxiv.org/abs/q-alg/9705011

\bibitem [P-T-1] {P-T-1} 
J.~H.~Przytycki, P.~Traczyk, Invariants of links of Conway type, 
{\em Kobe J. Math. } 4, 1987, 115-139.

\bibitem [P-Ts] {P-Ts}
J.~H.~Przytycki, T. Tsukamoto, The fourth skein module and 
the Montesinos-Nakanishi conjecture for 3-algebraic links, 
{\it  J. Knot Theory Ramifications}, 10(7), 2001, 959--982.\\
e-print:\ http://arxiv.org/abs/math.GT/0010282

\bibitem [Ron-2] {Ron-2}
Y. Rong, Link polynomials of higher order, {\it J.~London Math. Soc. 2}, 
56(1), 1997, 189--208.


\bibitem [R-S] {R-S}
C.~Rourke, B.~Sanderson, 
{\em Introduction to piecewise-linear topology}, 
Ergebnisse der Mathematik und ihrer Grenzgebiete, Band 69;
Springer-Verlag, New York-Heidelberg, 1972.

\bibitem [Si-3] {Si-3}
A.~S. Sikora, Skein modules and TQFT, Knots in Hellas' 98;
The Proceedings of the International Conference
on Knot Theory and its Ramifications; Volume 1.
In the Series on Knots and Everything, Vol. 24, September 2000.

\bibitem [Si-4] {Si-4}
A.~S. Sikora, $PSL_n$-character varieties as spaces of graphs,
{\it Trans. Amer. Math. Soc.}, 353, 2001, 2773-2804. \\
e-print:\ http://arxiv.org/abs/math.RT/9806016

\bibitem [Si-5] {Si-5}
A.~S. Sikora, Skein theory for SU(n)-quantum invariants;\\
http://front.math.ucdavis.edu/math.QA/0407299

\bibitem [T-L]{T-L}
H.~N. V.Temperley. E.~H. Lieb. Relations between the ``percolation" 
and ``coloring" problem and other graph-theoretical problems associated 
with regular planar lattices: 
some exact results for the ``percolation" problem. {\em Proc. Roy. Soc. 
Lond.} A 322, 1971, 251-280.
 
\bibitem[Tu-2]{Tu-2}
 V.~G.~Turaev, The Conway and Kauffman modules of the solid torus,
{\it Zap. Nauchn. Sem. Lomi} 167 (1988), 79-89. English translation:
{\it J. Soviet Math.} 52, 1990, 2799-2805.

\bibitem[Tu-3]{Tu-3} 
 V.~G.~Turaev, Skein quantization of Poisson algebras of loops 
on surfaces, {\em Ann. Scient. \'Ec. Norm. Sup.}, 4(24), 1991, 635-704.

\bibitem[Tu-4]{Tu-4}
V.~G. Turaev, Quantum invariants of knots and 3-manifolds,
de Gruyter Studies in Mathematics, 18,
Walter de Gruyter \& Co., Berlin, 1994. x+588 pp.

\bibitem[Tut-1]{Tut-1}
W.~T. Tutte. A ring in graph theory. {\it Proc. Cambridge Phil. Soc.},
 43, 1947, 26-40.

\bibitem[Wa-1]{Wa-1}
B. Wajnryb, A simple presentation for the mapping
class group of an orientable surface, {\it Israel  J. Math.},
45(2/3), 1983, 157-174.

\bibitem[Wa-2]{Wa-2}
B. Wajnryb, Markov classes in certain finite symplectic
representations of braid groups, in:  Braids (Santa Cruz, CA, 1986),
687--695, Contemp. Math., 78, Amer. Math. Soc.,
Providence, RI, 1988.

\bibitem[Wa-3]{Wa-3}
B. Wajnryb, A braidlike presentation of ${\rm Sp}(n,p)$,
{\it Israel J. Math.} 76(3), 1991, 265--288.

\bibitem[Wit]{Wit} E. Witten. Quantum field theory and the Jones polynomial.
{\em Comm. Math. Phys. } 121, 1989, 351-399.
\end{thebibliography}
\end{document}